\documentclass[12pt]{article}
\usepackage{latexsym,amssymb,amsmath,mathrsfs}
\usepackage{sectsty}
\usepackage[anchorcolor=blue,%
 bookmarks=true,%
 bookmarksnumbered=true,%
 dvipdfm,%
]{hyperref}

\allsectionsfont{\normalsize\sc\center}

\makeatletter \@addtoreset{equation}{section} \makeatother

\newcommand{\nn}{\nonumber}
\newcommand{\ds}{\displaystyle}

\newtheorem{thm}{Theorem}[section]
\newtheorem{prop}[thm]{Proposition}
\newtheorem{cor}[thm]{Corollary}
\newtheorem{lem}[thm]{Lemma}
\newtheorem{remark}[thm]{Remark}
\newtheorem{assumption}{Assumption}

\newenvironment{demo}[1]{\par\begin{trivlist}%
\item[]{\bf #1}\ }{\end{trivlist}\par}
\newcommand{\Proof}{\begin{demo}{\bf Proof.\ }}
\newcommand{\Proofof}[1]{\begin{demo}{\bf Proof of #1.\ }}
\newcommand{\bpf}{\Proof}
\newcommand{\apf}{\Proofof}
\newcommand{\epf}{\hfill $\square$ \end{demo}}

\newcommand{\Thm}[1]{Theorem~\ref{th:#1}}
\newcommand{\Prop}[1]{Proposition~\ref{prop:#1}}
\newcommand{\Lem}[1]{Lemma~\ref{lem:#1}}
\newcommand{\Cor}[1]{Corollary~\ref{cor:#1}}
\newcommand{\Rem}[1]{Remark~\ref{rem:#1}}
\newcommand{\Ass}[1]{Assumption~\ref{ass:#1}}
\newcommand{\eq}[1]{\eqref{eq:#1}}

\renewcommand{\P}{\mathbb{P}}
\newcommand{\E}{\mathbb{E}}
\newcommand{\R}{\mathbb{R}}
\newcommand{\N}{\mathbb{N}}
\newcommand{\sC}{\mathscr{C}}
\newcommand{\sD}{\mathscr{D}}
\newcommand{\sF}{\mathscr{F}}
\newcommand{\sL}{\mathscr{L}}

\renewcommand{\a}{\alpha}
\renewcommand{\b}{\beta}
\newcommand{\gm}{\gamma}
\newcommand{\dl}{\delta}
\newcommand{\ep}{\varepsilon}
\newcommand{\zt}{\zeta}
\newcommand{\h}{\eta}
\newcommand{\lm}{\lambda}
\newcommand{\kp}{\kappa}
\newcommand{\ro}{\rho}
\newcommand{\sg}{\sigma}
\newcommand{\ph}{\varphi}

\newcommand{\Dl}{\Delta}
\newcommand{\Lm}{\Lambda}
\newcommand{\Ph}{\Phi}

\newcommand{\abs}[1]{\left| #1 \right|}

\newcommand{\abra}[1]{\left( #1 \right)}
\newcommand{\bbra}[1]{\left\{ #1 \right\}}
\newcommand{\cbra}[1]{\left[ #1 \right]}
\newcommand{\dbra}[1]{\langle #1 \rangle}
\newcommand{\ebra}[1]{\lfloor #1 \rfloor}

\newcommand{\wg}{\wedge}
\newcommand{\nab}{\nabla}
\newcommand{\e}{\mathrm{e}}

\title
{\bf \Large
Coupling by reflection of diffusion processes 
via discrete approximation under a backward Ricci flow
\footnotetext{
    Partially supported 
    by the JSPS fellowship for research abroad. 
}
}

\author{
  Kazumasa Kuwada 
}

\DeclareMathOperator{\Ric}{Ric}

\DeclareMathOperator{\Cut}{Cut}
\newcommand{\Cutst}{\Cut_{\mathrm{ST}}}

\begin{document}
\maketitle
\begin{abstract}
A coupling by reflection of a time-inhomogeneous 
diffusion process on a manifold are studied. 
The condition we assume is a natural time-inhomogeneous 
extension of lower Ricci curvature bounds. 
In particular, 
it includes the case of backward Ricci flow. 
As in time-homogeneous cases, 
our coupling provides a gradient estimate of 
the diffusion semigroup 
which yields the strong Feller property. 
To construct the coupling via discrete approximation, 
we establish the convergence in law of geodesic random walks 
as well as a uniform non-explosion type estimate. 
\end{abstract}
\section{Introduction} 
\label{sec:intro} 


In stochastic analysis, 
coupling methods of stochastic processes 
have played a prominent role in the literature. 
Given two stochastic processes $Y_1 (t)$ and $Y_2 (t)$ 
on a state space $M$, 
a coupling $\mathbf{X} (t) = ( X_1 (t) , X_2 (t) )$ 
of $Y_1 (t)$ and $Y_2 (t)$ 
is a stochastic process on $M \times M$ such that 
$X_i$ has the same law as $Y_i$ for $i=1,2$. 
By constructing a suitable coupling which reflects 
the geometry of the underlying structure, 
one can obtain various estimates 
for heat kernels, harmonic maps, eigenvalues 
etc.~under natural geometric assumptions 
(see \cite{Hsu,Kend_survey,Wang_book05} for instance).
Recently, the heat equation 
on time-inhomogeneous spaces such as Ricci flow 
have been studied intensively 
(see 
\cite{Arn-Coul-Thal_horiz,
  Coul_gtBM,
  K-Phili, 
  K-Phili2, 
  McC-Topp_Wass-RF,
  Topp_Lopt,
  ZhangQS_HK-RF-Poinc} 
and references therein). 
These studies have succeeded in revealing 
a tighter connection 
between the heat equation and 
the underlying geometric structure even 
in time-inhomogeneous cases. 
It should be remarked that 
an idea of coupling methods lies behind some of them 
\cite{Arn-Coul-Thal_horiz,K-Phili2,McC-Topp_Wass-RF,Topp_Lopt} 
in an essential way . 


This paper is aimed at constructing 
a coupling by reflection of a diffusion process 
associated with a time-dependent family of metrics 
such as (backward) Ricci flow. 
Let $M$ be a smooth manifold 
with a family of complete Riemannian metrics 
$\{ g(t) \}_{t \in [ T_1 , T_2 ] }$. 
By $\{ X (t) \}_{t \in [ T_1 ,T_2 ] }$, 
we denote the $g(t)$-Brownian motion. 
It means that 
$X (t)$ is the time-inhomogeneous 
diffusion process on $M$ 
associated with $\Delta_{g(t)} /2$, 
where $\Delta_{g(t)}$ is the Laplacian 
with respect to $g(t)$ 
(see \cite{Coul_gtBM} for a construction of $g(t)$-Brownian motion). 
As in time-homogeneous cases studied in 
\cite{Crans,
  Kend,
  K8,
  Renes_poly,
  Wang94,
  Wang_book05} 
under a lower Ricci curvature bound, 
a coupling by reflection 
$\mathbf{X} (t) = ( X_1 (t) , X_2 (t) )$ of 
two $g(t)$-Brownian motions 
starting from a different point 
provides us a useful control of 
the coupling time $\tau^*$, 
the first time when coupled particles meet. 
A simple version of our main theorem 
which states such a control is as follows: 
\begin{thm} \label{th:main} 
Suppose 
\begin{equation} \label{eq:c-b} 
\partial_t g (t) 
\le 
\Ric_{g(t)} 
\end{equation} 
holds. 
Then, for each $x_1 , x_2 \in M$, 
there exists a coupling 
$\mathbf{X} (t) := ( X_1 (t) , X_2 (t) )$ 
of two $g(t)$-Brownian motions 
starting at $(x_1 , x_2 )$ 
satisfying 
\begin{equation} \label{eq:dom}
\P [ \tau^* > t ] 
\le 
\P \cbra{ 
  \inf_{T_1 \le s \le t} B (s) 
  > - \frac{d_{g(T_1)} ( x_1 , x_2 )}{2} 
} 
\end{equation} 
for each $t$, where $d_{g(T_1)}$ is the distance 
function on $M$ with respect to $g(T_1)$ and 
$B (t)$ is a 1-dimensional standard Brownian motion 
starting at the time $T_1$. 
\end{thm}
For the complete statement of our main theorem, 
see \Thm{main0}. 
There we also study a diffusion process 
which generalizes the $g(t)$-Brownian motion. 
The condition \eq{c-b} can be interpreted as 
a time-inhomogeneous analogue of nonnegative Ricci curvature
(see \Rem{O-U}). 
This condition is essentially the same as 
backward super Ricci flow in \cite{McC-Topp_Wass-RF} 
(Our condition is slightly different in constant 
since our $g(t)$-Brownian motion 
and hence the heat equation 
corresponds to $\Dl_{g(t)} / 2$ 
instead of $\Dl_{g(t)}$). 
Obviously, \eq{c-b} is satisfied if $g(t)$ evolves 
according to backward Ricci flow 
$\partial_t g(t) = \Ric_{g(t)}$. 
As in the time-homogeneous case, 
our coupling time estimate yields  
a gradient estimate of the heat semigroup 
which implies the strong Feller property 
for the heat semigroup (see \Cor{sF}). 
Note that, when $g(t)$ is a backward Ricci flow, 
the same estimate as \Cor{sF} is also obtained 
in \cite{Coul_gtBM} 
by using techniques in stochastic differential geometry. 

To explain our approach to \Thm{main}, 
let us review a heuristic idea of 
the construction of a coupling by reflection 
as well as that of the derivation of \eq{dom}. 
Given a Brownian particle $X_1$, 
we will construct $X_2$ by determining 
its infinitesimal motion 
$d X_2 (t) \in T_{X_2 (t)} M$ 
by using $d X_1 (t) \in T_{X_1 (t)} M$. 
First we take a minimal $g(t)$-geodesic 
$\gm$ joining $X_1 (t)$ and $X_2 (t)$. 
Next, by using the parallel transport along $\gm$ 
associated with the $g(t)$-Levi-Civita connection, 
we bring $d X_1 (t)$ into $T_{X_2 (t)} M$. 
Finally we define $d X_2 (t)$ 
as a reflection of it with respect to a hyperplane 
being $g(t)$-perpendicular to $\dot{\gm}$ 
in $T_{X_2 (t)} M$. 
From this construction, 
the It\^{o} formula implies that 
$d_{g(t)} ( X_1 (t) , X_2 (t) )$ 
should become a semimartingale at least until 
$( X_1 (t) , X_2 (t) )$ hits the $g(t)$-cutlocus $\Cut_{g(t)}$. 
The semimartingale decomposition is 
given by variational formulas of arclength. 
On the bounded variation part, 
there appear 
the time-derivative of $d_{g(t)}$ 
and the second variation of $d_{g(t)}$, 
which is dominated in terms of Ricci curvature. 
With the aid of our condition \eq{c-b}, 
these two terms are compensated and 
a nice domination of 
the bounded variation part follows. 
Thus the hitting time to 0 of $d_{g(t)} ( X_1 (t) , X_2 (t) )$, 
which is the same as $\tau^*$, can be 
estimated by that of the dominating semimartingale. 
Indeed, we can regard $2 B (t) + d_{g(T_1)} (x_1 , x_2 )$ 
which appeared in the right hand side of \eq{dom} 
as the dominating semimartingale. 
The effect of our reflection appears 
in the martingale part $2 B (t)$ 
which makes it possible 
for the dominating martingale to hit $0$. 
This construction seems to work 
as long as $( X_1 (t) , X_2 (t) )$ is not in the cutlocus. 
Moreover, we can hope it possible 
to construct it beyond the cutlocus 
since the cutlocus is sufficiently small and 
the effect of singularity at the cutlocus 
should make $d_{g(t)} ( X_1 (t) , X_2 (t) )$ to decrease. 
If we succeed in doing so, the bounded variation part 
will involve a ``local time at $\Cut_{g(t)}$''. 
It will be negligible since it would be nonpositive. 
We can conclude that 
almost all technical difficulties are 
concentrated on the treatment of singularity at the cutlocus 
in order to make this heuristic argument rigorous. 
In fact, 
\Thm{main} is shown in \cite{Phili_sem} 
by using SDE methods 
when the $g(t)$-cutlocus is empty 
for every $t \in [ T_1 , T_2 ]$. 

Our construction of a coupling by reflection 
is based on a time-discretized approximation 
as studied in \cite{K8,Renes_poly}. 
We construct a coupling of geodesic random walks 
each of whose marginals approximates 
the original diffusion process. 
The construction will be finished after taking a limit 
so that these approximations converge. 
Our method has a remarkable advantage 
in treating singularities arising from the cutlocus. 
In our construction, 
we can avoid to extract a local time at $\Cut_{g(t)}$ 
and directly obtain a dominating process 
which does not involve such a term. 
In the present framework, the singular set $\Cut_{g(t)}$ 
also depends on time parameter $t$ and hence treating it 
by using stochastic differential equations 
seems to be more complicated than 
in the time-homogeneous case. 

Different kinds of couplings are studied 
in above-mentioned papers. 
Based on the theory of optimal transportation, 
McCann and Topping 
\cite{McC-Topp_Wass-RF,Topp_Lopt} 
studied a coupling of heat distributions 
which minimizes their transportation cost. 
They used the squared distance in \cite{McC-Topp_Wass-RF} 
or 
Perelman's $\mathcal{L}$-functional in \cite{Topp_Lopt} 
respectively 
to quantify a transportation cost. 
Their coupling is 
closely related to coupling of Brownian motions 
by parallel transport 
along minimal ($\mathcal{L}$-)geodesics. 
In fact, 
studying a coupling by parallel transport 
by probabilistic methods 
recovered and extended (a part of) their results 
in \cite{Arn-Coul-Thal_horiz} and 
\cite{K-Phili2} respectively. 
Note that 
our approach via time-discretized approximation 
is used in \cite{K-Phili2}. 
In addition, 
we also can construct a coupling by parallel transport 
by using our method to recover 
a result in \cite{Arn-Coul-Thal_horiz} 
(see \Thm{parallel}). 
It explains that our approach is also effective 
even when we study a different kind of couplings. 

We give a remark on 
a difference in methods between ours and 
Arnaudon, Coulibaly and Thalmaier's one 
\cite{Arn-Coul-Thal_horiz} 
to construct a coupling by parallel transport. 
They consider one-parameter family of 
coupled particles along a curve. 
Intuitively saying, 
they concatenate coupled particles along a curve 
by iteration of making a coupling by parallel transport. 
Since ``adjacent'' particles are 
infinitesimally close to each other, 
we can ignore singularities on the cutlocus 
when we construct a coupled particle 
from an ``adjacent'' one. 
It should be noted that 
their method does not seem to be able to 
be applied directly 
in order to construct a coupling by reflection. 
Indeed, their construction of a chain of coupled particles 
heavily relies on a multiplicative (or semigroup) property 
of the parallel transport. 
However, our reflection operation obviously 
fails to possess such a multiplicative property.  
Since our reflection map changes orientation,  
there is no chance to interpolate it 
with a continuous family of isometries. 

In what follows, we will state 
the organization of this paper. 
In the next section, 
we show basic properties 
of a family of Riemannian manifolds 
$( ( M , g(t) ) )_{t}$. 
In particular, we prove that 
Riemannian metrics $( g(t) )_t$ are 
locally comparable with each other. 
It will be used to give a uniform control of 
several error terms which appear 
as a result of our discrete approximation. 
In section 3, we will study geodesic random walks 
in our time-inhomogeneous framework. 
There we introduce them and 
prove the convergence in law 
to a diffusion process. 
After a small discussion 
at the beginning of the section, 
the proof is divided into two main parts. 
In the first part, 
we will give a uniform estimate 
for the exit time from a big compact set of 
geodesic random walks. 
Our assumption here is almost the same as 
in \cite{K-Phili} where non-explosion 
of the diffusion process is studied 
(see \Rem{tight} \ref{KP} for more details). 
In the second part, we prove tightness of geodesic random walks 
on the basis of the result in the first part. 
In section 4, we will construct a coupling by reflection 
and show an estimate of coupling time, 
which completes the proof of \Thm{main} as a special case. 


%
\section{Properties on time-dependent metric} 
\label{sec:pre} 
Let $M$ be a $m$-dimensional manifold. 
For $-\infty < T_1 < T_2 < \infty$, 
Let $( g(t) )_{t \in [ T_1 , T_2 ]}$ 
be a family of complete Riemannian metrics on $M$ 
which smoothly depends on $t$. 
\begin{remark} \label{rem:time-parameter} 
It seems to be restrictive that 
our time parameter only runs over 
the compact interval $[ T_1 , T_2 ]$. 
An example of $g(t)$ we have in mind is 
a solution to the backward Ricci flow equation. 
In this case, we \emph{can} work 
on a semi-infinite interval $[ T_1 ,\infty )$ 
only when we study an ancient solution of the Ricci flow. 
Thus $T_2 < \infty$ is not so restrictive. 
In addition, we could extend 
our results to the case on $[ T_1 , \infty )$ 
with a small modification of our arguments. 
It would be helpful to study an ancient solution. 
To deal with a singularity of Ricci flow, 
it could be nice to work 
on a semi-open interval $( T_1 , T_2 ]$, 
where $T_1$ is the first time 
when a singularity emerges. 
In that case, we should be more careful 
since we cannot give 
``an initial condition at $T_1$'' 
to define a $g(t)$-Brownian motion on $M$. 
\end{remark}

We collect some notations 
which will be used in the sequel. 
Throughout this paper, 
we fix a reference point $o \in M$. 
Let $\N_0$ be nonnegative integers. 
For $a, b \in \R$, 
$a \wg b$ and $a \vee b$ stand 
for $\min\{ a, b \}$ and $\max\{ a , b \}$ 
respectively. 
Let $\Cut_{g(t)} (x)$ be 
the set of the $g(t)$-cutlocus of $x$ on $M$. 
Similarly, the $g(t)$-cutlocus $\Cut_{g(t)}$ and 
the space-time cutlocus $\Cutst$ are 
defined by 
\begin{align*}
\Cut_{g(t)} 
& : = 
\bbra{ 
  (x,y) \in M \times M 
  \; | \; 
  y \in \Cut_{g(t)} (x)
},
\\ 
\Cutst 
& : = 
\bbra{ 
  (t,x,y) \in [ T_1 , T_2 ] \times M \times M 
  \; | \; 
  (x,y) \in \Cut_{g(t)} 
}.
\end{align*} 
Set $D(M) := \{ (x , x) \; | \; x \in M \}$. 
The distance function with respect to $g(t)$ 
is denoted by $d_{g(t)} (x,y)$.  
Note that 
$\Cutst$ is closed 
and that $d_{g(\cdot)} ( \cdot, \cdot )$ is 
smooth on 
$
[ T_1 , T_2 ] \times M \times M 
\setminus 
( \Cutst \cup [ T_1 , T_2 ] \times D(M) )
$
(see \cite{McC-Topp_Wass-RF}, cf. \cite{K-Phili}). 
We denote an open $g(s)$-ball of radius $R$ 
centered at $x \in M$ 
by $B^{(s)}_R (x)$. 
Some additional notations will be given 
at the beginning of the next section.

In the following three lemmas 
(\Lem{metric_control}-\Lem{bdd_cpt}), 
we discuss a local comparison between 
$d_{g(t)}$ and $d_{g(s)}$ for $s \neq t$.  
Those will be a geometric basis of 
the further arguments. 
%
%
\begin{lem} \label{lem:metric_control}
Let $M_0$ be a compact subset of $M$. 
Then there exists $\kp = \kp ( M_0 )$ 
such that 
\[
\e^{ - 2 \kp | t-s | } g (s) 
\le 
g (t) 
\le 
\e^{ 2 \kp | t-s | } g (s) 
\] 
holds on $M_0$ for $t,s \in [ T_1 , T_2 ]$. 
In particular, 
if 
a minimal $g(s)$-geodesic $\gm$ joining $x,y \in M_0$ 
is included in $M_0$, 
then, for $t \in [ T_1 , T_2 ]$, 
\[
d_{g(t)} (x,y) \le \e^{\kp | t-s |} d_{g(s)} (x,y). 
\] 
\end{lem}
\bpf
Let $\pi \: : \: TM \to M$ be a canonical projection. 
Let us define $\hat{M}_0$ by 
\[
\hat{M}_0 
: =
\bbra{ 
  ( t, v ) \in [ T_1 , T_2 ] \times TM
  \; \left| \; 
      \pi (v) \in M_0 , 
      \abs{ v }_{g(t)} \le 1 
  \right.
}. 
\]
Note that $\hat{M}_0$ is closed 
since $g(\cdot)$ is continuous. 
We claim that $\hat{M}_0$ is sequentially compact. 
Let us take a sequence 
$( ( t_n , v_n ) )_{n \in \N} \subset \hat{M}_0$. 
We may assume $t_n \to t \in [ T_1 , T_2 ]$ 
and $\pi (v_n) \to p \in M_0$ as $n \to \infty$ 
by taking a subsequence if necessary. 
Let $U$ be a neighborhood of $p$ such that 
$\{ v \in TM \; | \; \pi (v ) \in U \} 
\simeq 
U \times \R^m$. 
For sufficiently large $n$, 
we regard $v_n$ as an element of $U \times \R^m$ 
and 
write $v_n = ( p_n , \tilde{v}_n )$. 
If we cannot take any convergent subsequence of 
$( v_n )_{n \in \N}$, 
then $| \tilde{v}_n | \to \infty$ 
as $n \to \infty$, 
where $| \cdot |$ stands for 
the standard Euclidean norm on $\R^m$ 
(irrelevant to $( g(t) )_{t \in [ T_1 , T_2 ]}$). 
Set $v_n' = ( p_n , | \tilde{v}_n |^{-1} \tilde{v}_n )$. 
Then, there exists a subsequence 
$( v_{n_k}' )_{k \in \N} \subset ( v_n' )_{n \in \N}$ 
such that 
$v_{n_k}' \to v_{\infty}' = (p, \bar{v}' )$ 
as $n \to \infty$ 
for some $\bar{v}' \in \R^m$ with $| \bar{v}' | =1$. 
Since $g(\cdot)$ is continuous, 
$
g( t_{n_k} ) ( v_{n_k}' , v_{n_k}' ) 
\to 
g(t)( v_{\infty}' , v_{\infty}' )
$
as $k \to \infty$. 
On the other hand, 
$
g ( t_{n_k} ) ( v_{n_k}' , v_{n_k}' ) 
\le 
| \tilde{v}_{n_k} |^{-2} \to 0
$
since $g (t_n ) ( v_n , v_n ) \le 1$. 
Thus $\bar{v}'$ must be 0. 
It contradicts with $| \bar{v}'| = 1$. 
Hence $\hat{M}_0$ is sequentially compact. 

Since 
$
\hat{M}_0 \ni ( t, v ) 
\mapsto 
\partial_t g (t) (v,v)
$ 
is continuous, 
there exists a constant $\kp = \kp ( M_0 ) >0$
such that 
$
\abs{ 
\partial_t g (t) ( v, v )
} 
\le 2 \kp 
$ 
for every $( t , v ) \in \hat{M}_0$. 
Take $v \in \pi^{-1}( M_0 )$, $v \neq 0_{\pi (v)}$. 
Then 
\[ 
\partial_t g (t) ( v, v ) 
= 
\abs{v}_{g(t)}^2 
\partial_t g (t) 
( 
 \abs{v}_{g(t)}^{-1} v , 
 \abs{v}_{g(t)}^{-1} v 
) 
\le 
2 \kp \abs{v}_{g(t)}^2
. 
\]
Thus $\partial_t \log g (t) (v,v) \le 2 \kp$ holds.  
By integrating it from $s$ to $t$ with $s < t$, 
we obtain 
$g(t) (v,v) \le \e^{2 \kp ( t-s )} g (s) (v,v)$. 
We can obtain the other inequality similarly. 

For the latter assertion, 
for $a,b$ with $\gm (a) = x$ and $\gm (b) = y$,  
\[
d_{g(t)} (x,y) 
\le 
\int_a^b \abs{\dot{\gm}(u)}_{g(t)} du 
\le 
\e^{\kp \abs{ t-s } } 
\int_a^b \abs{\dot{\gm}(u)}_{g(s)} du 
= 
\e^{\kp |t-s|} d_{g(s)} (x,y)
. 
\]
\epf
\begin{lem} \label{lem:ball_control} 
For $R > 0$, $x \in M$ and $t \in [ T_1 , T_2 ]$, 
there exists 
$\dl = \dl ( x, t, R ) > 0$ such that 
$\bar{B}^{(s)}_{r} (x) \subset \bar{B}^{(t)}_{3r} (x)$ 
for $r \le R$ and $s \in [ T_1 , T_2 ]$ 
with $\abs{ s - t } \le \dl$. 
\end{lem}
\bpf
Set $\kp := \kp ( \bar{B}^{(t)}_{3R} (x) )$ 
as in \Lem{metric_control} and 
$\dl : = \kp^{-1} \log 2$. 
Take $p \in \bar{B}^{(s)}_{r} (x)$ and 
a minimal $g(s)$-geodesic 
$\gm \: : \: [a,b] \to M$ 
joining $x$ and $p$. 
Suppose that there exists $u_0 \in [a,b]$ 
such that $\gm (u_0) \in \bar{B}^{(t)}_{3r} (x)^c$. 
Let 
$
\bar{u}_0 
: = 
\inf 
\{ 
  u \in [ a,b ] 
  \; | \; 
  \gm ( u ) \in \bar{B}^{(t)}_{3r} (x)^c 
\}
$. 
Since 
$
\gm ( [ a , \bar{u}_0 ] ) 
 \subset 
\bar{B}^{(t)}_{3r} (x) 
 \subset 
\bar{B}^{(t)}_{3R} (x)
$ 
and $d_{g(t)} ( x , \gm ( \bar{u}_0 ) ) = 3r$, 
\Lem{metric_control} yields 
\[
d_{g(s)} (x,p) 
\ge 
\int_a^{\bar{u}_0} \abs{\dot{\gm}(u)}_{g(s)} d u
\ge 
\e^{-\kp \dl}
\int_a^{\bar{u}_0} \abs{\dot{\gm}(u)}_{g(t)} d u 
= 
\frac{3r}{2}. 
\]
This is absurd. 
Hence $\gm ([a,b]) \in \bar{B}^{(t)}_{3r} (x)$. 
In particular, 
$\gm (b) = p \in \bar{B}^{(t)}_{3r} (x)$. 
\epf
\begin{lem} \label{lem:bdd_cpt}
For $R > 0$, 
there exists a compact subset 
$M_0 = M_0 (R)$ of $M$ 
such that 
\begin{equation} \label{eq:bdd_cpt}
\bbra{ 
  p \in M 
  \; \left| \;  
      \inf_{t \in [ T_1 , T_2 ] } 
      d_{g(t)} ( o , p ) \le R 
  \right.
} 
\subset M_0 . 
\end{equation}
\end{lem}
\bpf
For each $t \in [ T_1 , T_2 ]$, 
take $\dl ( o, t, R + 1 ) > 0$ 
according to \Lem{ball_control}. 
Take 
$\{ t_i \}_{i=1}^n \subset [ T_1 , T_2 ]$ 
such that 
\[ 
[ T_1 , T_2 ] 
\subset 
\bigcup_{i=1}^n 
( 
 t_i - \dl ( o, t_i , R + 1) , 
 t_i + \dl ( o, t_i , R + 1) 
) 
.
\] 
Let us define a compact set $M_0 \subset M$ 
by 
$
M_0 := 
\bigcup_{i=1}^n \bar{B}^{(t_i)}_{3R} (o) 
$. 
Take $p \in M$ such that  
$\inf_{T_1 \le t \le T_2} d_{g(t)} (o,p) \le R$. 
For $\ep \in (0,1)$, 
take $s \in [ T_1 , T_2 ]$ 
such that $d_{g(s)} (o,p) \le R + \ep$.
Then there exists $j \in \{ 1, \ldots, N \}$ 
such that $|s - t_j | < \dl ( o, t_j , R + 1)$. 
By \Lem{ball_control}, 
it implies 
$
p 
\in \bar{B}^{(s)}_{R + \ep} (o) 
 \subset 
\bar{B}^{(t_j)}_{3( R + \ep )} (o) 
 \subset 
\bigcup_{i=1}^n 
\bar{B}^{(t_i)}_{3( R + \ep )} (o) 
$. 
Hence the conclusion follows 
by letting $\ep \downarrow 0$. 
\epf 
Another useful consequence of 
\Lem{metric_control} and \Lem{ball_control} 
is the following: 
\begin{lem} \label{lem:d-conti}
$d_{g(\cdot)} (\cdot, \cdot)$ is continuous 
on $[T_1 , T_2 ] \times M \times M$. 
\end{lem}
\bpf
Since the topology 
on $[T_1 ,T_2 ] \times M \times M$ 
is metrizable, 
It suffices to show 
$
\lim_{n\to \infty} d_{g(t_n)} ( x_n , y_n ) 
 = 
d_{g(t)} (x,y)
$ when $(t_n , x_n , y_n ) \to ( t , x , y )$ 
as $n \to \infty$. 
By the triangle inequality, 
\begin{equation} \label{eq:d-conv} 
\abs{ 
  d_{g(t_n)} ( x_n , y_n ) 
  - 
  d_{g(t)} ( x , y ) 
} 
\le 
\abs{ 
  d_{g(t_n)} (x,y) - d_{g(t)} (x,y) 
} 
+ d_{g(t_n)} (x,x_n) 
+ d_{g(t_n)} (y,y_n) 
. 
\end{equation}
Take $R > 0$ so that $B_{R}^{(t)} (x)$ includes 
a minimal $g(t)$-geodesic joining $x$ and $y$. 
Take $\kp = \kp ( \bar{B}_{4R}^{(t)} (x) )$ 
according to \Lem{metric_control}. 
We can easily see that 
every minimal $g(t)$-geodesic 
joining $y$ and $y_n$ is included in $B_{2R}^{(t)} (x)$ 
for sufficiently large $n \in \N$. 
Thus \Lem{metric_control} yields 
\[
\limsup_{n\to\infty} d_{g(t_n)} (y,y_n) 
\le 
\limsup_{n \to \infty} 
\e^{\kp | t - t_n|} d_{g(t)} (y,y_n) 
= 0.   
\] 
We can show $d_{g(t_n)} (x,x_n) \to 0$ similarly. 
Take a minimal $g(t_n)$-geodesic 
$\gm_n \: : \: [ a, b ] \to M$ 
joining $x$ and $y$. 
By our choice of $R$, 
\Lem{metric_control} again yields 
\[
d_{g(t_n)} ( x, \gm_n (u) ) 
 \le 
d_{g(t_n)} ( x, y ) 
 \le 
\e^{\kp | t - t_n | } d_{g(t)} (x,y) 
 \le 
\e^{\kp | t - t_n | } R. 
\]
It implies 
$\limsup_{n\to\infty} d_{g(t_n)} (x,y) \le d_{g(t)} (x,y)$. 
In addition, 
$\gm_n$ is included in $B_{4R/3}^{(t_n)} (x)$ 
for sufficiently large $n$. 
Thus \Lem{ball_control} and \Lem{metric_control} 
yield 
$d_{g(t)} (x,y) \le \e^{\kp | t - t_n |} d_{g(t_n)} (x,y)$. 
Hence the conclusion follows 
by combining these estimates with \eq{d-conv}. 
\epf 
Before closing this section, 
we will provide a local lower bound 
of injectivity radius 
which is uniform in time parameter. 
\begin{lem} \label{lem:inject}
For every $M_1 \subset M$ compact, 
there is $\tilde{r}_0 = \tilde{r}_0 ( M_1 ) > 0$ 
such that 
$d_{g(t)} (y,z) < \tilde{r}_0$ implies 
$(t,y,z) \notin \Cutst$ 
for any 
$(t,y,z) \in [ T_1 , T_2 ] \times M_1 \times M_1$. 
\end{lem}
\bpf 
Take $R > 1$ so that 
$
\sup_{ t \in [ T_1 , T_2 ] } 
\sup_{ x \in M_1 } 
d_{g(t)} ( o , x ) 
< R - 1
$. 
By \Lem{bdd_cpt}, 
there exists a compact set $M_0 \subset M$ 
such that \eq{bdd_cpt} holds. 
For every $t \in [ T_1 , T_2 ]$ and 
$x \in M_1$, 
$(t,x,x) \notin \Cutst$. 
It implies that 
there is $\h_{t,x} \in ( 0, 1 )$ 
such that 
$( s, y, z ) \notin \Cutst$ 
whenever 
$
d_{g(t)} (x,y) 
\vee 
d_{g(t)} (x,z) 
\vee 
|t-s| 
< \h_{t,x}
$ 
since $\Cutst$ is closed. 
Thus there exist $N \in \N$ and 
$( t_i , x_i ) \in [ T_1 , T_2 ] \times M_1$ 
($i = 1, \ldots , N$)   
such that 
\[
[ T_1 , T_2 ] \times M_1  
\subset 
\bigcup_{i=1}^N 
\abra{ 
  t_i - \frac{\h_{t_i , x_i}}{2} , 
  t_i + \frac{\h_{t_i , x_i}}{2} 
} 
\times 
B_{\h_{t_i , x_i} / 2 }^{(t_i)} (x_i) 
.
\] 
Set $\tilde{r}_0 > 0$ by 
\[
\tilde{r}_0 := 
\frac12 
\exp 
\abra{ 
  - \frac{\kp}{2} \max_{1 \le i \le N} \h_{t_i , x_i}  
} 
\min_{1 \le i \le N} \h_{t_i , x_i} 
, 
\]
where $\kp = \kp ( M_0 ) > 0$ is 
as in \Lem{metric_control}. 
Take $(s,y,z) \in [ T_1 , T_2 ] \times M_1 \times M_1$ 
with $d_{g(s)} ( y, z ) < \tilde{r}_0$. 
Take $j \in \{ 1 , \ldots, N \}$ 
so that 
$
| s - t_j | 
\vee  
d_{g(t_j)} ( x_j , y )  
< \h_{t_j , x_j} / 2 
$. 
By virtue of 
the choice of $R$ and $M_0$, 
\Lem{bdd_cpt} yields that 
every $g(s)$-geodesic joining $y$ and $z$ 
is included in $M_0$. 
Thus \Lem{metric_control} yields 
\[
d_{g(t_j)} ( y , z )
\le 
\e^{\kp | s - t_j |} d_{g(s)} ( y , z )
< 
\frac{\h_{t_j , x_j}}{2}.  
\]
It implies 
$
|s - t_j| 
\vee d_{g(t_j)} ( x_j , y ) 
\vee d_{g(t_j)} ( x_j , z )
< \h_{t_j , x_j}
$ and hence $(s,y,z) \notin \Cutst$. 
\epf
\section{Approximation via geodesic random walks} 
\label{sec:IP}

Let $( Z(t) )_{t \in [ T_1 , T_2 ]}$ be 
a family of smooth vector fields 
continuously depending 
on the parameter $t \in [ T_1 , T_2 ]$. 
Let $X(t)$ be the diffusion process 
associated with the time-dependent generator 
$\sL_t = \Dl_{g(t)}/2 + Z(t)$ 
(see \cite{Coul_gtBM} for a construction of $X(t)$ 
by solving a SDE on the frame bundle). 
Note that $( t , X(t) )$ is a unique solution 
to the martingale problem associated with 
$\partial_t + \sL_\cdot$ 
on $[ T_1 ,T_2 ] \times M$ 
(see \cite{Hsu} for the time-homogeneous case. 
Its extension to time-inhomogeneous case is straightforward; 
see \cite{Str-Var} also). 

In what follows, 
we will use several notions in Riemannian geometry 
such as exponential map $\exp$, Levi-Civita connection $\nab$, 
Ricci curvature $\Ric$ etc. 
To clarify the dependency on the metric $g(t)$, 
we put $(t)$ on superscript or $g(t)$ on subscript. 
For instance, we use the following symbols: 
$\exp^{(t)}$, $\nab^{(t)}$ and $\Ric_{g(t)}$. 
We refer to \cite{Chavel2} 
for basics in Riemannian geometry 
which will be used in this paper. 

For each $t \in [ T_1 , T_2 ]$, 
we fix a measurable section 
$\Ph^{(t)} \: : \: M \to \mathscr{O}^{(t)} (M)$ 
of the $g(t)$-orthonormal frame bundle 
$\mathscr{O}^{(t)} (M)$ of $M$. 
Take a sequence of independent, identically distributed 
random variables $\{ \xi_n\}_{n\in \N}$ which are 
uniformly distributed on the unit disk in $\R^m$. 
Given $x_0 \in M$, 
let us define 
a continuously-interpolated 
geodesic random walk 
$( X^\a ( t ) )_{t \in [ T_1 , T_2 ]}$ 
on $M$ 
starting from $x_0$ 
with a scale parameter $\a > 0$ 
inductively. 
Let $t_n^{(\a)} : = ( T_1 + \a^2 n ) \wg T_2$ 
for $n \in \N_0$. 
For $t = T_1 = t_0^{(\a)}$, 
set $X^\a (T_1) := x_0$. 
after $X^\a (t)$ is defined 
for $t \in [ T_1 , t_n^{(\a)} ]$,  
we extend it to $t \in [ t_n^{(\a)} , t_{n+1}^{(\a)} ]$ 
by 
\begin{align} 
\nn 
\tilde{\xi}_{n+1} 
& : = 
\sqrt{m+2} \Ph^{( t_n^{(\a)} )} ( X^\a ( t_n^{(\a)} ) ) \xi_{n+1}
, 
\\ \nn 
X^\a ( t ) 
& := 
\exp_{X^\a ( t_n^{(\a)} )}^{( t_n^{(\a)} )} 
\bigg(
  \frac{ t - t_n^{(\a)} }{\a^2}  
  \Big( 
    \a \tilde{\xi}_{n+1} 
     + 
    \a^2 Z( t_n^{(\a)} ) 
  \Big)
\bigg) . 
\end{align}
%
For later use, 
we define 
$
N^{(\a)} 
: = 
\inf 
\bbra{
  n \in \N_0 
  \; | \; 
  t_{n+1} - t_n < \a^2 
}
$. 
This is the total number of discrete steps of 
our geodesic random walks 
with scale parameter $\a$. 
Set $\sC : = C ([ T_1 , T_2 ] \to M)$ 
and $\sD : = D ( [ T_1 , T_2 ] \to M)$. 
By using a distance $d_{g(T_1)}$ on $M$, 
we metrize $\sC$ and $\sD$ as usual 
so that $\sC$ and $\sD$ become Polish spaces 
(see \cite{Ethier-Kurtz} 
for a distance function on $\sD$, for example). 
Set $\sC_1 := C ( [ T_1 , T_2 ] \to [ 0 , \infty ) )$. 
%
Let us define a time-dependent 
$(0,2)$-tensor field $( \nab Z (t))^{\flat}$  
by 
\[
( \nab Z (t))^\flat (X,Y) 
: = 
\frac12 
\abra{ 
  \dbra{ \nab_X^{(t)} Z (t) , X }_{g(t)} 
  + 
  \dbra{ \nab_Y^{(t)} Z(t) , Y }_{g(t)}
}
.  
\]
\begin{assumption} \label{ass:non-explosion} 
There exists a locally bounded 
nonnegative measurable function $b$ 
on $[ 0 , \infty )$ such that 
\begin{enumerate}
\item
For all $t \in [ T_1 , T_2 )$, 
\[
2 ( \nabla Z (t) )^\flat +  \partial_t g (t) 
\le 
\Ric_{g(t)} + b ( d_{g(t)} ( o , \cdot ) ) g(t)
.
\] 
\item \label{non-explosion}
For each $C > 0$, 
a 1-dimensional diffusion process $y_t$ given by 
\[
d y_t 
= 
d \beta_t 
+ 
\frac12 
\left( 
    C 
    + 
    \int_0^{y_t} b ( s ) ds 
\right) 
dt
\] 
does not explode. 
(This is the case if and only if
\[
\int_1^\infty 
\exp \abra{ 
    - \int_1^y \mathbf{b} (z) dz 
} 
\int_1^y 
\exp \abra{ 
  \int_1^z \mathbf{b} (\xi) d \xi 
} dz 
dy 
= \infty, 
\]
where $\mathbf{b}(y) := C + \int_0^{y} b(s) ds$. 
see e.g.~\cite[Theorem~VI.3.2]{Ik-Wat}.)
\end{enumerate}
\end{assumption} 
Our goal in this section is to prove the following: 
\begin{thm} \label{th:IP}
Under \Ass{non-explosion}, 
$X^\a$ converges in law to $X$ in $\sC$ as $\a \to 0$.  
\end{thm}
Most of arguments in this section will be 
devoted to show the tightness i.e. 
\begin{prop} \label{prop:tight} 
$( X^\a )_{\a \in (0,1)}$ is tight in $\sC$. 
\end{prop}
In fact, as we will see in the following, 
\Prop{tight} easily implies \Thm{IP}. 
\apf{\Thm{IP}}
By virtue of \Prop{tight}, 
for any subsequence of $( X^\a )_{\a \in (0,1)}$
there exists a further subsequence $ ( X^{\a_k} )_{k \in \N}$ 
which converges in law in $\sC$ as $k \to \infty$. 
Thus it suffices to show that 
this limit has the same law as $X$. 
Let $( \b^\a (t) )_{t \in [0,\infty)}$ be 
a Poisson process of intensity $\a^{-2}$ 
which is independent of $\{ \xi_n \}_{n \in \N}$. 
Set 
\[
\bar{\b}^\a (t) 
 := 
( T_1 + \a^2 \b^\a ( t - T_1 ) ) 
 \wg 
t_{N^{(\a)}}^{(\a)}
. 
\]
Then the Poisson subordination 
$X^{\a_k} ( \bar{\b}^{\a_k} (\cdot) )$ also 
converges in law in $\sD$ to the same limit 
(see \cite{Billi} for instance). 
Note that 
$( 
 \bar{\b}^{\a} (t) , 
 X^{\a} ( \bar{\b}^{\a} (t) ) 
)_{t \in [ T_1 , T_2 ]}$ 
is a space-time Markov process. 
The associated semigroup 
$P^{(\a)}_t$ and its generator $\tilde{\sL}^{(\a)}$
are given by 
\begin{align*}
P^{(\a)}_t f 
& := 
\e^{ - ( t - T_1 )\a^{-2}} 
\abra{ 
  \sum_{l=1}^{ N^{(\a)} }
  \frac{ ( ( t - T_1 ) \a^{-2} )^l }{l!} 
  (q^{(\a)} )^l f  
  + 
  \sum_{l >  N^{(\a)} } 
  \frac{ ( ( t - T_1 ) \a^{-2} )^l }{l!} 
  ( q^{(\a)} )^{ N^{(\a)} } f 
}, 
\\
\tilde{\sL}^{(\a)} f 
& := 
\a^{-2} ( q^{(\a)}f - f) , 
\intertext{where}
q^{(\a)} f (t,x)  
& := 
\E \cbra{ 
  f ( 
  t + \a^2 , 
  \exp^{(t)}_x 
  \abra{ 
    \a \sqrt{m+2} 
    \Ph^{(t)}(x) \xi_1 
    + 
    \a^2 Z ( t ) 
  }
}. 
\end{align*}
We can easily prove 
$\tilde{\sL}^{(\a)} f \to ( \partial_t + \sL_\cdot ) f$ 
uniformly as $\a \to 0$ 
for $f \in C^\infty_0 ( [ T_1 , T_2 ] \times M )$. 
Since 
$( 
 \bar{\b}^{\a} (t) , 
 X^{\a} (  \bar{\b}^{\a} (t) ) 
)_{t \in [ T_1 , T_2 ]}$ 
is a solution to the martingale problem 
associated with $\tilde{\sL}^{(\a)}$, 
the limit in law of 
$( 
 \bar{\b}^{\a_k} (t) , 
 X^{\a} (  \bar{\b}^{\a_k} (t) ) 
)_{t \in [ T_1 , T_2 ]}$ 
solves the martingale problem associated 
with $\partial_t + \sL_\cdot$. 
By the uniqueness of the martingale problem, 
this limit has the same law as 
that of $( t , X(t) )_{t \in [ T_1 , T_2 ]}$. 
It completes the proof.
\epf

\begin{remark} \label{rem:tight}
\begin{enumerate}
\item
A result on a convergence of semigroups 
\cite{Kurtz_appr-semigr} 
was used to show the convergence of 
finite dimensional distributions 
in the time-homogeneous case 
\cite{Blum} (see \cite{Renes_poly} also). 
It is not so clear that we can employ the same argument 
in our time-inhomogeneous case. 
One difficulty arises from the absence of invariant measures 
for semigroups even in the case $Z(t) \equiv 0$. 
Although the $g(t)$-Riemannian measure is 
a unique invariant measure for $\Dl_{g(t)}$, 
this measure also depends on time parameter. 
Thus we cannot expect that 
it becomes an invariant measure of semigroups. 
This obstacle also prevents us to employ the existing theory of 
time-dependent Dirichlet forms 
(see \cite{Oshima_tdepDF} for instance) 
in order to study our problem. 
\item \label{KP}
\Prop{tight} also asserts that 
any subsequential limit in law 
is a probability measure on $\sC$. 
Since we have not added any cemetery point to $M$ 
in the definition of $\sC$, 
\Thm{IP} implies that $X$ cannot explode. 
It almost recovers the result in \cite{K-Phili}. 
Our assumption is slightly stronger than that in \cite{K-Phili} 
on the point where we require \ref{non-explosion} 
for \emph{all} $C > 0$, not a given constant. 
Note that 
we will use \Ass{non-explosion} \ref{non-explosion} 
only for a specified constant $2 C_0$ given in \Lem{drift}. 
However, its expression looks complicated and 
it seems to be less interesting 
to provide a explicit bound. 
\end{enumerate}
\end{remark}
Now we introduce some additional notations 
which will be used in the rest of this paper. 
For $t \in [ T_1 , T_2 ]$, 
we define $\ebra{ t }_\a$ by 
$\ebra{ t }_\a 
:= 
\sup \bbra{ 
  \a^2 n + T_1 
  \; | \; 
  n \in \N_0 , \a^2 n + T_1 < t 
}$. 
Set $\sF_n : = \sg ( \xi_1 , \ldots \xi_n )$. 
For $R > 1$, 
let us define 
$\sg_R \: : \: \sC_1 \to [ T_1 , T_2 ] \cup \{ \infty \}$ 
by 
\[
\sg_R (w)
: = 
\inf 
\bbra{ 
  t \in [ T_1 , T_2 ] 
  \; | \; 
  w (t) > R - 1
},
\] 
where $\inf \emptyset = \infty$. 
We write  
$\hat{\sg}_R 
:= 
\sg_R ( d_{g (\cdot)} ( o , X^\a (\cdot) ) ) 
$ and 
$
\bar{\sg}_R 
: = 
\a^{-2} ( \ebra{ \hat{\sg}_R }_\a - T_1 ) + 1
$. 
Note that $\bar{\sg}_R$ is an $\sF_n$-stopping time.  
For each $t \in [ T_1 , T_2 ]$ 
and $x , y \in M$ with $x \neq y$, 
we choose 
a minimal unit-speed $g(t)$-geodesic 
$\gm_{xy}^{(t)} \: : \: [ 0, d_{g(t)} (x,y) ] \to M$ 
from $x$ to $y$. 
Note that 
we can choose $\gm_{xy}^{(t)}$ so that 
$(x,y) \mapsto \gm_{xy}^{(t)}$ is measurable 
in an appropriate sense 
(see e.g. \cite{Renes_poly}). 
We use the same symbol $\gm_{xy}^{(t)}$ 
for its range $\gm_{xy}^{(t)} ( [ 0, d_{g(t)} (x,y)] )$. 
\subsection{Uniform bound for escape probability}
\label{sec:Uni-prob}

The goal of this subsection is to show the 
following: 
\begin{prop} \label{prop:non-explosion}
$ 
\lim_{R \uparrow \infty} \limsup_{\a \downarrow 0} 
\P \cbra{ \hat{\sg}_R \le T_2 } = 0
$. 
\end{prop} 
For the proof, 
we will establish a discrete analogue of 
a comparison argument 
for the radial process 
as discussed in \cite{K-Phili}. 
In this subsection, 
we fix $R > 1$ sufficiently large 
so that $d_{g(T_1)} ( o , x_0 ) < R-1$ 
until the final line of the proof of 
\Prop{non-explosion}. 
We also fix a relatively compact open set 
$M_0 \subset M$ satisfying \eq{bdd_cpt}. 
Set $r_0 := \tilde{r}_0 \wg (1/2)$, 
where $\tilde{r}_0 = \tilde{r}_0 ( M_0 )$ 
is as in \Lem{inject}. 

The first step for proving \Prop{non-explosion} 
is to show 
a difference inequality for the radial process 
$d_{g(t)} ( o , X^\a (t) )$ (\Lem{variation}). 
It will play a role of 
the It\^{o} formula for the radial process 
in our discrete setting. 
We introduce some notations 
to discuss how to avoid the singularity of 
$d_{g(\cdot)} ( o , \cdot )$ 
on $\{ o \} \cup \Cut_{g(\cdot)} (o)$. 
For $r > 0$, 
let us define a set $A_r' , A_r''$ and $A_r$ 
as follows: 
\begin{align*} 
A_r' 
& : = 
\bbra{ 
  (t,x,y) \in [ T_1 , T_2 ] \times M_0 \times M_0 
  \; \left| \; 
      \begin{array}{l}
          d_{g(t)} (x,x') + d_{g(t)} (y,y') + \abs{ t- t' }
          \ge r 
          \\
          \mbox{for any $(t', x', y') \in \Cutst$}
      \end{array}
  \right.
}, 
\\
A_r'' 
& : = 
\bbra{ 
  (t,x,y) \in [ T_1 , T_2 ] \times M_0 \times M_0 
  \; | \; 
  d_{g(t)} (x,y) \ge r 
},  
\\
A_r 
& := 
A_r' \cap A_r'' 
.  
\end{align*} 
Note that $A_r$ is compact 
and that 
$d_{g(\cdot)} ( \cdot, \cdot )$ is smooth 
on $A_r$. 
For $t \in [ T_1 , T_2 ]$ and $p \in M$, 
let us define $o_p^{(t)} \in M_0$ 
by 
\[
o_p^{(t)} 
: = 
\begin{cases}
\ds 
\gm_{o p}^{(t)} 
\abra{ \frac{r_0}{2} } 
& 
\mbox{
  if $( t , o , p ) \notin A_{r_0}'$, 
} 
\\
o & \mbox{otherwise.} 
\end{cases}
\]
For simplicity of notations, 
we denote $o_{X^\a ( t_n^{(\a)} )}^{( t_n^{(\a)} )}$ by $o_n$. 
Similarly, we use the symbol $\gm_n$ for 
$\gm_{o_n X^\a ( t_n^{(\a)} )}^{( t_n^{(\a)} )}$ 
throughout this section. 
Note that 
$( t , o_{p}^{(t)} , p ) \notin \Cutst$ 
holds. 
Furthermore, 
it is uniformly separated from $\Cutst$ 
in the following sense: 
%
\begin{lem} \label{lem:away0}
There exist $r_1 > 0$ and $\dl_1 > 0$ 
such that the following holds: 
Let $t_0 , t \in [ T_1 , T_2 ]$ 
with $t - t_0 \in [ 0 , \dl_1]$. 
Let $p_0 \in B_{R-1}^{(t_0)} (o)$ and 
$p \in B_{\dl_1}^{(t_0)} (p_0)$. 
Then we have 
\begin{enumerate}
\item \label{approx_error0}
$
d_{g(t)} ( o , p ) 
\le 
\e^{\kp (t - t_0)} 
\abra{ d_{g(t_0)} ( o , p_0 ) + d_{g(t_0)} ( p_0 , p ) }
$, 
\item \label{away_cutlocus0}
$( t , o_{p_0}^{(t_0)} , p ) \in A_{r_1}$ 
when $p_0 \notin B_{r_0}^{(t_0)} (o)$. 
\end{enumerate}
Here $\kp = \kp (  M_0 ) > 0$ is given 
according to \Lem{metric_control}.
\end{lem}
By applying \Lem{away0} to $X^\a$, 
we obtain the following: 
\begin{cor} \label{cor:away}
There exist $\a_0 > 0$ and 
$h \: : \: [ 0 , \a_0 ] \to [ 0, 1 ]$ 
with $\lim_{\a \downarrow 0} h (\a) = 0$ 
such that the following holds: 
For $\a \le \a_0$, $n \in \N_0$ and 
$s,t \in [ t_n^{(\a)}, t_{n+1}^{(\a)} ]$, 
when $n < \bar{\sg}_R$, 
\begin{enumerate}
\item \label{approx_error}
$
d_{g(t)} ( o , X^\a (s) ) 
\le 
\e^{ \kp \a^2 } 
\abra{ 
  d_{g( t_n^{(\a)} )} ( o , X^\a ( t_n^{(\a)} ) ) 
  + 
  h (\a) 
}
$, 
\item 
$( t , o_n , X^\a ( s ) ) \in A_{r_1}$ 
when $X^\a ( t_n^{(\a)} ) \notin B^{(t_n^{(\a)})}_{r_0} (o)$. 
\end{enumerate}
Here $r_1$ is the same as in \Lem{away0}. 
\end{cor}
\apf{\Cor{away}}
Set 
$
\bar{Z} 
: = 
\sup_{t \in [ T_1 , T_2 ], x \in M_0} 
\abs{ Z(t) }_{g(t)} (x) 
$. 
Note that we have 
$
d_{g( t_n^{(\a)} )} ( X^\a ( t_n^{(\a)} ) , X^\a (t) ) 
\le 
\sqrt{m+2} \a + \bar{Z} \a^2
$ 
by the definition of $X^\a$.  
Take $\a_0 > 0$ so that 
$\sqrt{m+2} \a_0 + \bar{Z} \a_0^2 \le \dl_1$ 
and 
$\a^2 \le \dl_1$ 
hold, where $\dl_1$ is as in \Lem{away0}. 
Then the conclusion follows 
by applying \Lem{away0} 
with $t_0 = t_n^{(\a)}$, 
$p_0 = X^\a ( t_n^{(\a)} )$ and 
$p = X^\a ( s )$. 
\epf
%
%
\apf{\Lem{away0}}
We show that 
\ref{approx_error0} holds with $\dl_1 = 1$. 
By the triangle inequality, 
the proof is reduced to showing 
the following two inequalities: 
\begin{align} \label{eq:d-1}
d_{g(t)} ( o , p_0 ) 
& \le 
\e^{\kp ( t - t_0 )} 
d_{g( t_0 )} ( o , p_0 )
,
\\ \label{eq:d-2}
d_{g(t)} ( p_0 , p ) 
& \le 
\e^{\kp (t - t_0)} 
d_{g( t_0 )} ( p_0, p )
.
\end{align}
Our condition \eq{bdd_cpt} yields that 
$\gm_{op_0}^{(t_0)}$ is included in $M_0$. 
Thus \Lem{metric_control} yields \eq{d-1}. 
When $p \in B_1^{(t_0)} (p_0)$, 
we have $p \in B_R^{(t_0)} (o)$.  
Hence 
\eq{bdd_cpt} and \Lem{metric_control} 
yield \eq{d-2} 
in a similar way as \eq{d-1}. 

Let us turn to consider \ref{away_cutlocus0}. 
For simplicity of notations, 
we denote $o_{p_0}^{(t_0)}$ by $o'$ 
in this proof. 
We assume 
that $t - t_0 \in [ 0 , \dl ]$ and 
$p \in B_\dl^{(t_0)} (p_0)$ hold 
for $\dl > 0$. 
First we will show 
$(t , o' , p ) \in A_{r_0 / 4}''$ 
when $\dl$ is sufficiently small. 
Note that 
$( t_0 , o' , p_0 ) \in A_{r_0 / 2}''$ holds 
since $p_0 \notin B^{(t_0)}_{r_0} (o)$ and 
$d_{g( t_0 )} ( o , o' ) \in \{ r_0 / 2 , 0 \}$. 
Let $q \in \gm_{o' p_0}^{(t)}$. 
By the triangle inequality, 
\begin{equation} \label{eq:d-3}
d_{g(t)} ( o , q ) 
\le 
d_{g(t)} ( o , o' )
+  
d_{g(t)} ( o' , p_0 ) 
.
\end{equation}
Since $r_0 /2 < 1 < R$ holds, 
\eq{bdd_cpt} yields 
$\gm^{(t_0)}_{oo'} \subset M_0$ 
when $o' \neq o$. 
We can easily see that 
$\gm_{o'p_0}^{(t_0)} \subset \gm_{o p_0}^{(t_0)} \subset M_0$.  
Thus, 
by applying \Lem{metric_control} to \eq{d-3}, 
\begin{align} \nn
d_{g(t)} ( o , q ) 
& \le  
\e^{\kp (t - t_0)} 
\abra{ 
  d_{g(t_0)} ( o , o' ) 
  + 
  d_{g(t_0)} ( o' , p_0 )
} 
\\ \label{eq:d-4}
& = (R-1) \e^{\kp \a^2}. 
\end{align}
Take $\dl_2 := 1 \wg ( \kp^{-1} \log (R / (R-1) ) )$. 
Then, for any $\dl \in ( 0 , \dl_2 )$, 
\eq{d-4} and \eq{bdd_cpt} imply 
$\gm_{o'p_0}^{(t)} \subset M_0$. 
Hence the triangle inequality, 
\Lem{metric_control} and 
\eq{d-2} yield 
\begin{align} \nn
d_{g(t)} ( o' , p ) 
& \ge 
d_{g(t)} ( o' , p_0 ) 
-
d_{g(t)} ( p_0 , p )
\\ \nn
& \ge 
\e^{-\kp (t - t_0)} 
d_{g(t_0)} ( o' , p_0 ) 
- 
\e^{\kp (t - t_0)} 
d_{g(t_0)} ( p_0 , p )
\\ \label{eq:d-5}
& \ge 
\frac{\e^{-\kp \dl} r_0}{2} 
- 
\e^{\kp \dl} \dl 
\end{align} 
when $\dl \le \dl_2$. 
Thus there exists 
$
\dl_3 = 
\dl_3 ( \kp , r_0 , R)
\in 
(0, \dl_2 ]
$ 
such that 
the right hand side of \eq{d-5}
is greater than $r_0 / 4$ 
whenever $\dl \in ( 0 , \dl_3 )$. 
Hence 
$( t , o' , p ) \in A_{r_0 / 4}''$ 
holds in such a case. 

Next we will show 
that there exists $r_1' > 0$ 
such that 
$( t , o' , p ) \in A_{r_1'}'$ 
holds for sufficiently small $\dl$.  
Once we have shown it, 
the conclusion holds with 
$r_1 = r_1' \wg ( r_0 / 4 )$. 
As we did in showing 
$(t , o' , p ) \in A''_{r_0 / 4}$, 
we begin with studying 
the corresponding statement 
for $( t_0 , o' , p_0 )$. 
More precisely, we claim that 
there exists 
$r_1'' \in ( 0, 1 )$ such that 
$ 
( t_0 , o' , p_0 ) 
\in A_{r_1''} 
$ 
for all $\dl \in ( 0, 1 )$. 
When $o' = o$, 
$( t_0 , o' , p_0 ) \in A_{r_0}'$ 
directly follows 
from the definition of $o' = o_{p_0}^{(t_0)}$. 
When $o' \neq o$, 
set 
\[ 
H : = 
\bbra{ 
  (t,x,y) \in [ T_1 , T_2 ] \times \bar{M}_0 \times \bar{M}_0 
  \; \left| \; 
      \begin{array}{l}
          r_0 \le d_{g(t)} (o,y) \le R-1 , \, 
          d_{g(t)} (o,x) = r_0 / 2 ,  
          \\
          d_{g(t)} (x,y) 
          = 
          d_{g(t)} (o,y) - d_{g(t)} (o,x)
      \end{array} 
  \right.
} . 
\]
Note that $H$ is compact and 
that $H \cap \Cutst = \emptyset$ holds 
since $(t,x,y) \in H$ implies that 
$x$ is on a minimal $g(t)$-geodesic 
from $y$ to $o$. 
Since 
$( t_0 , o' , p_0 ) \in H$ 
by the definition of $o'$, 
it suffices to show that 
there exists $\tilde{r}_1 > 0$ 
such that $H \subset A_{\tilde{r}_1}'$. 
Indeed, the claim will be shown 
with $r_1'' = \tilde{r}_1 \wg r_0$ 
once we have proved it. 
Suppose that 
$H \subset A_r'$ does not hold 
for any $r \in ( 0 , 1 )$. 
Then there are sequences 
$(t_j , x_j , y_j ) \in H$, 
$(t_j' , x_j' , y_j' ) \in \Cutst$, 
$j \in \N$ such that 
$
| t_j - t_j' | 
+ d_{g(t_j)} (x_j , x_j') 
+ d_{g(t_j)} (y_j , y_j')
\to 0$ 
as $j \to \infty$. 
We may assume 
that $( (t_j , x_j , y_j ) )_j$ converges. 
Since $(t_j , x_j , y_j ) \in H$, 
$x_j' , y_j' \in M_0$ holds 
for sufficiently large $j$. 
Thus we can take a convergent subsequence of 
$( (t_j' , x_j' , y_j' ) )_j$. 
Since 
$\Cutst$ and $H$ are closed and 
$d_{g(\cdot)} ( \cdot ,  \cdot )$ is continuous, 
it contradicts with $H \cap \Cutst = \emptyset$. 
 
To complete the proof, we show that 
there exists $\dl_1 \in ( 0 , \dl_3 ]$ 
such that 
$( t , o' , p ) \in A_{r_1'' / 2}'$ 
when $\dl \in ( 0 , \dl_1 )$. 
Suppose that there exists 
$(t' , x' , y' ) \in \Cutst$ 
such that 
$
| t - t' | 
+ d_{g(t)} ( o' , x' ) 
+ d_{g(t)} ( p , y' ) 
< r_1'' / 2
$. 
For any $q \in \gm_{py'}^{(t)}$, 
the triangle inequality and 
the assertion \ref{approx_error0} yield 
\begin{equation} \label{eq:d-6} 
d_{g(t)} ( o , q ) 
\le 
d_{g(t)} ( o , p ) 
+ 
d_{g(t)} ( p , y' ) 
\le 
\e^{ \kp \dl } 
\abra{ 
  R - 1 
  + 
  \dl 
}
+ r_1'' / 2 
. 
\end{equation}
A similar observation implies 
$d_{g(t)} ( o , q' ) \le ( \e^{\kp \dl} r_0 + r_1'' ) / 2$ 
for $q' \in \gm_{o'x'}^{(t)}$. 
Thus there is $\dl_4 = \dl_4 ( \kp , R ) \in ( 0 , \dl_3 ]$ 
such that the right hand side of \eq{d-6}
is less than $R$ 
and $( \e^{\kp \dl} r_0 + r_1'' ) / 2 \le R$ 
whenever $\dl \in ( 0 , \dl_4 )$. 
In such a case, 
$\gm_{py'}^{(t)} \subset M_0$ 
and 
$\gm_{o'x'}^{(t)} \subset M_0$ 
hold. 
Since $( t_0 , o' , p_0 ) \in A_{r_1''}'$, 
\Lem{metric_control} yields 
\begin{align} \nn
| t - t' | 
+ & d_{g(t)} ( o' , x' ) 
+ d_{g(t)} ( p , y' )
\\ \nn
& \ge 
| t_0 - t' | - \dl 
+ \e^{-\kp \dl} d_{g(t_0)} ( o' , x' ) 
+ \e^{-\kp \dl} d_{g(t_0)} ( p , y' )
\\ \label{eq:d-7}
& \ge 
\e^{- \kp \dl} r_1'' 
+ ( 1 - \e^{- \kp \dl} ) | t_0 - t' | 
- \dl 
- \e^{- \kp \dl} \dl 
\end{align}
Take $\dl_1 = \dl_1 ( \kp , r_1'' ) \in ( 0 , \dl_4 ]$ 
so that the right hand side of \eq{d-7} 
is greater than $r''_1 / 2$ when $\dl \in ( 0 , \dl_1 )$. 
Then \eq{d-7} is absurd 
for any $\dl \in ( 0 , \dl_1 )$. 
Thus it implies the conclusion. 
\epf

We prepare some notations 
for the second variational formula for the arclength. 
Let $\nabla^{(t)}$ be the $g(t)$-Levi-Civita connection 
and $\mathcal{R}^{(t)}$ the $g(t)$-curvature tensor 
associated with $\nabla^{(t)}$. 
For a smooth curve $\gm$ and 
smooth vector fields $U, V$ along $\gm$, 
the index form $I_\gm^{(t)} ( U , V )$ is given by 
\[
I_\gm^{(t)} ( U , V )
: = 
\int_\gm 
\abra{
  \dbra{ 
    \nab_{\dot{\gm}}^{(t)} U , 
    \nab_{\dot{\gm}}^{(t)} V 
  }_{g(t)} 
  - 
  \dbra{ 
    \mathcal{R}^{(t)} ( U, \dot{\gm} ) \dot{\gm}, V 
  }_{g(t)}
} ds . 
\]
We write 
$I_\gm^{(t)} (U,U) =: I_\gm^{(t)} (U)$ 
for simplicity of notations. 
Let $G_{t,x,y} (u)$ be the solution to 
the following initial value problem 
on $[ 0, d (x,y) ]$: 
\begin{equation} \nn 
\begin{cases}
\ds 
G_{t,x,y}'' ( u ) 
= 
- 
\frac{ 
  \Ric_{g(t)} 
  ( \dot{\gm}_{xy}^{(t)} (u) , \dot{\gm}_{xy}^{(t)} (u) ) 
}{
  m-1
}
G_{t,x,y} (u)
, 
\\
\ds 
G_{t,x,y} (0) = 0 , \, G_{t,x,y}' (0) = 1 . 
\end{cases}
\end{equation}
Note that $G_{t,x,y} (u) > 0$ for $u \in ( 0, d(x,y) ]$ 
if $y \notin \Cut_{g(t)} (x)$ 
(see \cite[Proof of Lemma~9]{K-Phili}). 
For simplicity, 
we write $G_n : = G_{t_n^{(\a)} , o_n , X^\a ( t_n^{(\a)} )}$. 
When $X^\a ( t_n^{(\a)} ) \notin B_{r_0}^{( t_n^{(\a)} )} (o)$, 
we define a vector field $V^\dag$ along $\gm_n$ 
for each $V \in T_{X^\a ( t_n^{(\a)} )} M$ 
by 
\[ 
V^\dag (\gm_n (u))
: = 
\frac{ 
  G_n (u) 
}{
  G_n ( d_{g( t_n^{(\a)} )} ( o_n , X^\a ( t_n^{(\a)} ) ) )
}
( /\!\!/_{\gm_n}^{(t_n^{(\a)})} V )( \gm_n (u) ),
\] 
where $/\!\!/_{\gm_n}^{(t_n^{(\a)})} V$ 
is the parallel vector field  
along $\gm_n$ of $V$ 
associated with $\nab^{( t_n^{(\a)} )}$. 
Take $v \in \R^m$. 
By using these notations, 
for $n \in \N_0$ with $n < N^{(\a)}$, 
let us define $\lm_{n+1}$ and $\Lm_{n+1}$ 
by 
\begin{align*}
\lm_{n+1} 
& : = 
\dbra{ 
  \tilde{\xi}_{n+1} 
  , 
  \dot{\gm}_n 
}_{g ( t_n^{(\a)} )} 
, 
\\
\Lm_{n+1} 
& : = 
\partial_t d_{g( t_n^{(\a)} )} ( o , o_n ) 
+ 
\partial_t d_{g( t_n^{(\a)} )} ( o_n , X^\a ( t_n^{(\a)} ) ) 
\\
& \qquad + 
\dbra{ 
  Z ( t_n^{(\a)} ) 
  , 
  \dot{\gm}_n 
}_{g ( t_n^{(\a)} )} ( X^\a ( t_n^{(\a)} ) ) 
+ 
\frac{1}{2} 
I_{\gm_n}^{( t_n^{(\a)} )} 
( 
  \tilde{\xi}_{n+1}^\dag 
)
\end{align*}
when 
$
X^\a ( t_n^{(\a)} ) 
\notin 
B^{( t_n^{(\a)} )}_{r_0} (o) 
$, 
and 
$\lm_{n+1} = \sqrt{m+2} \dbra{ \xi_{n+1} , v }_{\R^m}$ 
and $\Lm_{n+1} = 0$ otherwise. 
%
%
\begin{lem} \label{lem:variation}
If $n < \bar{\sg}_R \wg N^{(\a)}$, 
$\a < \a_0$ is sufficiently small 
and 
$X^\a ( t_n^{(\a)} ) \notin B_{r_0}^{( t_n^{(\a)} )} (o)$, 
then 
\begin{equation*}
d_{g( t_{n+1}^{(\a)} )} ( o, X^\a ( t_{n+1}^{(\a)} ) ) 
\le 
d_{ g ( t_n^{(\a)} ) } ( o , X^\a ( t_n^{(\a)} ) ) 
+ \a \lm_{n+1} 
+ \a^2 \Lm_{n+1} 
+ o ( \a^2 ) 
\end{equation*}
almost surely, 
where $\a_0$ is as in \Cor{away}. 
In addition, 
$o( \a^2 )$ is controlled uniformly. 
\end{lem}
\bpf
By virtue of \Cor{away}, 
for sufficiently small $\a$, 
the Taylor expansion together with 
the second variational formula 
yields 
\begin{align} \nn
d_{g( t_{n+1}^{(\a)} )} ( o_n , X^\a ( t_{n+1}^{(\a)} ) ) 
& \le 
d_{ g ( t_n^{(\a)} ) } ( o_n , X^\a ( t_n^{(\a)} ) ) 
+ \a \lm_{n+1} 
+ \a^2 \partial_t d_{g( t_n^{(\a)} )} ( o_n , X^\a ( t_n^{(\a)} ) ) 
\\ \nn 
& \quad 
+ 
\a^2 \dbra{ 
  Z ( t_n^{(\a)} ) 
  , 
  \dot{\gm}_n 
}_{g ( t_n^{(\a)} )} ( X^\a ( t_n^{(\a)} ) ) 
+ 
\frac{\a^2}{2} 
I_{\gm_n}^{( t_n^{(\a)} )} 
( 
  J_{\tilde{\xi}_{n+1}}
)
\\ \label{eq:2nd} 
& \quad 
+ o ( \a^2 ), 
\end{align} 
where $J_{\tilde{\xi}_{n+1}}$ is 
a $g( t_n^{(\a)} )$-Jacobi field 
along $\gm_n$ 
with the boundary value condition 
$J_{\tilde{\xi}_{n+1}} ( o_n ) = 0$ and 
$J_{\tilde{\xi}_{n+1}} ( X^\a ( t_n^{(\a)} ) ) = \tilde{\xi}_{n+1}$. 
Note that $o ( \a^2 )$ can be chosen uniformly 
since this expansion can be done 
on the compact set $A_{r_1}$ and 
every geodesic variation is included in $M_0$. 
By the index lemma, 
we have 
$
I_{\gm_n}^{( t_n^{(\a)} )} ( J_{\tilde{\xi}_{n+1}} ) 
\le 
I_{\gm_n}^{( t_n^{(\a)} )} ( \tilde{\xi}_{n+1}^\dag )
$. 
Hence the desired inequality follows when $o_n = o$. 
In the case $o_n \neq o$, we have 
\begin{align*}
d_{g( t_{n+1}^{(\a)} )} ( o , X^\a ( t_{n+1}^{(\a)} ) ) 
& \le 
d_{g( t_{n+1}^{(\a)} )} ( o , o_n ) 
+  
d_{g( t_{n+1}^{(\a)} )} ( o_n , X^\a ( t_{n+1}^{(\a)} ) ), 
\\
d_{g( t_n^{(\a)} )} ( o , X^\a ( t_n^{(\a)} ) ) 
& = 
d_{g ( t_n^{(\a)} )} ( o , o_n ) 
+ 
d_{g( t_n^{(\a)} )} ( o_n , X^\a ( t_n^{(\a)} ) ).
\end{align*}
Note that $(t_n^{(\a)} , o, o_n )$ is uniformly away from $\Cutst$ 
because of our choice of $r_0$ and \Lem{inject}. 
Therefore the conclusion follows 
by combining them with \eq{2nd}. 
\epf 
Before turning into the next step, 
we show the following two complementary lemmas 
(\Lem{LLN} and \Lem{drift}) 
which provide a nice control of 
the second order term $\Lm_n$ 
in \Lem{variation}. 
Set 
$\bar{\Lm}_n = \E [ \Lm_n \, | \, \sF_{n-1} ]$. 
%
\begin{lem} \label{lem:LLN} 
Let $( a_n )_{ n \in \N_0}$ be a uniformly bounded 
$\sF_n$-predictable process. 
Then 
\[
\lim_{\a \to 0}  
\a^2 
\sup 
\bbra{ 
  \left. 
  \Big| 
    \sum_{j=n}^{N+1} 
    a_j ( \Lm_j - \bar{\Lm}_j ) 
  \Big| 
  \; \right| \; 
  n , N \in \N, 
  n \le N \le N^{(\a)} \wg \bar{\sg}_R 
} 
= 0 
\quad 
\mbox{in probability.} 
\] 
\end{lem}
\bpf
Note that the map 
$(t,x,y) \mapsto G_{t,x,y} (d(x,y))$ 
is continuous on $A_{r_1}$. 
Since we have $G_{t,x,y} (d(x,y)) > 0$ on $A_{r_1}$, 
there exists $K > 0$ such that 
$K^{-1} < G_{t,x,y} (d(x,y)) < K$. 
This fact together with \Cor{away} yields   
$\abs{ \Lm_j }$ and $\abs{ \bar{\Lm}_j }$ 
are uniformly bounded 
if $j < \bar{\sg}_R$. 
Since 
$\sum_{j=1}^n a_j ( \Lm_j - \bar{\Lm}_j )$ is 
an $\sF_n$-martingale and 
$\bar{\sg}_R$ is $\sF_n$-stopping time, 
the Doob inequality yields 
\begin{equation} \label{eq:Doob}
\lim_{\a \to 0}  
\a^2 
\sup_{0 \le N \le N^{(\a)} \wg \bar{\sg}_R} 
\abs{ 
  \sum_{j=1}^{N+1} 
  a_j ( \Lm_j - \bar{\Lm}_j ) 
} 
= 0 
\quad 
\mbox{in probability.} 
\end{equation}
Here we used 
the fact $\lim_{\a \to 0} \a^2 N^{(\a)} = T_2 - T_1$. 
Note that  
\begin{align*}
\bigcup_{ N = 1 }^{N^{(\a)} \wg \bar{\sg}_R} 
\bigcup_{ n = 1 }^N 
&
\bbra{ 
  \a^2 
  \abs{ 
    \sum_{ j = n }^{N+1} 
    a_j ( \Lm_j - \bar{\Lm}_j ) 
  } 
  > \dl 
}
\\
& 
\subset 
\bigcup_{ N = 1 }^{N^{(\a)} \wg \bar{\sg}_R} 
\bigcup_{ n = 2 }^N 
\bbra{ 
  \a^2 
  \abs{ 
    \sum_{ j = 1 }^{n-1} 
    a_j ( \Lm_j - \bar{\Lm}_j ) 
  } 
  > \frac{\dl}{2}  
}
\cup
\bbra{ 
  \a^2 
  \abs{ 
    \sum_{ j = 1 }^{N+1} 
    a_j ( \Lm_j - \bar{\Lm}_j ) 
  } 
  > \frac{\dl}{2}  
}
\\
& 
= 
\bbra{ 
  \a^2 
  \sup_{0 \le N \le N^{(\a)} \wg \bar{\sg}_R} 
  \abs{ 
    \sum_{j=0}^N 
    a_j ( \Lm_j - \bar{\Lm}_j ) 
  } 
  > \frac{\dl}{2}  
}
. 
\end{align*}
Thus the conclusion follows from \eq{Doob}. 
\epf
\begin{lem} \label{lem:drift}
There exists a deterministic 
constant $C_0 > 0$ being independent of 
$\a$ and $R$ such that the following holds: 
\[ 
\bar{\Lm}_{n+1} 
\le 
C_0 
+ 
\frac12 
\int_0^{d_{g( t_n^{(\a)} )} (o , X^\a ( t_n^{(\a)} ) )} 
b ( u ) \, du 
. 
\] 
\end{lem}
\bpf
By using 
$(m+2) 
\E [ \dbra{ \xi_n , e_i } \dbra{ \xi_n , e_j } ] 
= 
\dl_{ij} 
$, 
we obtain 
\begin{align} \nn 
\E \cbra{ 
  I_{\gm_n}^{( t_n^{(\a)} )} ( \tilde{\xi}_{n+1}^\dag ) 
} 
& = 
\sum_{j=2}^m  
I_{\gm_n}^{( t_n^{(\a)} )} 
\abra{ 
  \abra{ 
    \Phi^{( t_n^{(\a)} )} ( X^\a ( t_n^{(\a)} ) ) e_j 
  }^\dag 
} 
\\ \nn 
& = 
\frac{ 
  (m-1) 
  G_n' ( d ( o_n , X^\a ( t_n^{(\a)} ) ) ) 
}{
  G_n ( d ( o_n , X^\a ( t_n^{(\a)} ) ) ) 
}
.
\end{align}
Note that we have 
\begin{align*}
\dbra{
  Z ( t_n^{(\a)} ) , \dot{\gm}_n 
}_{g ( t_n^{(\a)} )} ( X^\a ( t_n^{(\a)} )) 
& - 
\dbra{
  Z ( t_n^{(\a)} ) , \dot{\gm}_n 
}_{g ( t_n^{(\a)} )} ( o_n ) 
\\
& = 
\int_{0}^{d_{g( t_n^{(\a)} )} ( o_n , X^\a ( t_n^{(\a)} ))} 
\partial_s 
\dbra{
  Z ( t_n^{(\a)} ) , \dot{\gm}_n 
}_{g ( t_n^{(\a)} )} ( \gm_n (s) )|_{s = u} 
\, du  
\\
& = 
\int_{0}^{d_{g( t_n^{(\a)} )} ( o_n , X^\a ( t_n^{(\a)} ))} 
\dbra{
  \nab_{\dot{\gm}_n}^{( t_n^{(\a)} )} Z( t_n^{(\a)} ) 
  , 
  \dot{\gm}_n 
}_{g ( t_n^{(\a)} )} ( \gm_n (u) ) 
\, du . 
\end{align*}
Recall that, 
for 
$
( t , x , y ) 
\notin 
\Cutst 
$, 
we have  
\[
\partial_t d_{g(t)} ( x, y ) 
= 
\frac12 \int_0^{d_{g(t)} (x,y)} 
( \partial_t g(t) )
\abra{ 
  \dot{\gm}^{(t)}_{xy} (u) 
  , 
  \dot{\gm}^{(t)}_{xy}(u) 
} 
du 
\]
(cf. \cite[Remark~6]{McC-Topp_Wass-RF}). 
By combining them with \Ass{non-explosion}, 
\begin{align} \nn 
\bar{\Lm}_{n+1} 
& = 
\partial_t d_{g ( t_n^{(\a)} )} ( o , o_n ) 
+ 
\frac12 
\int_{0}^{d_{g( t_n^{(\a)} )} ( o_n , X^\a ( t_n^{(\a)} ))} 
\partial_t g ( t_n^{(\a)} ) 
( 
 \dot{\gm}_n (u) , 
 \dot{\gm}_n (u) 
) du 
\\ \nn
& \quad + 
\dbra{ Z ( t_n^{(\a)} ) , \dot{\gm}_n }_{g( t_n^{(\a)} )} 
( X^\a ( t_n^{(\a)} ) ) 
+ 
\frac{ 
  (m-1) 
  G_n' ( d ( o_n , X^\a ( t_n^{(\a)} ) ) ) 
}{
  2 G_n ( d ( o_n , X^\a ( t_n^{(\a)} ) ) ) 
}
\\ \nn
& \le 
\frac12 
\int_{0}^{d_{g( t_n^{(\a)} )} ( o , X^\a ( t_n^{(\a)} ))} 
b ( u ) \, du 
+ 
\partial_t d_{g( t_n^{(\a)} )} ( o , o_n ) 
+ 
\dbra{ Z ( t_n^{(\a)} ) , \dot{\gm}_n }_{g ( t_n^{(\a)} )} 
( o_n ) 
\\ \nn
& 
\quad + 
\frac12 
\int_{0}^{d_{g( t_n^{(\a)} )} ( o_n , X^\a ( t_n^{(\a)} ))} 
\Ric_{g ( t_n^{(\a)} )} 
( 
 \dot{\gm}_n (u) , 
 \dot{\gm}_n (u) 
) du 
\\ \label{eq:expectation}
& \quad + 
\frac{ 
  (m-1) 
  G_n' ( d ( o_n , X^\a ( t_n^{(\a)} ) ) ) 
}{
  2 G_n ( d ( o_n , X^\a ( t_n^{(\a)} ) ) ) 
}. 
\end{align}
Here we used the fact $ b (u) \ge 0$
in the case $o_n \neq o$. 
Note that 
\[
\int_{0}^{r} 
\Ric_{g ( t_n^{(\a)} )} 
( 
 \dot{\gm}_n (u) , 
 \dot{\gm}_n (u) 
) du 
+ 
\frac{ 
  (m-1) 
  G_n' ( r ) 
}{
  G_n ( r ) 
} 
\]
is non-increasing as a function of $r$. 
Indeed, we can easily verify it by taking 
a differentiation. 
Set 
\begin{align*}
C_1 
& : = 
\sup_{t \in [ T_1 , T_2 ]}
\sup_{ x \in B_{r_0} ^{(t)} (o) } 
\abra{ 
  \abs{ Z ( t ) }_{g(t)} (x)
  + 
  \sup_{
    \begin{subarray}{l} 
        V \in T_x M 
        \\
        \abs{ V }_{g(t)} \le 1 
    \end{subarray}
  }
  \abra{
    \partial_t g(t) ( V , V ) 
    + 
    | \Ric_{g(t)} ( V , V ) | 
  }
} 
. 
\end{align*} 
By virtue of \Lem{metric_control}, 
$C_1 < \infty$ holds. 
By applying a usual comparison argument 
to $G_n' ( r_0 ) / G_n ( r_0 )$, 
we obtain 
\begin{align*}
\int_{0}^{d_{g( t_n^{(\a)} )} ( o_n , X^\a ( t_n^{(\a)} ))} 
& \Ric_{g ( t_n^{(\a)} )} 
( 
 \dot{\gm}_n (u) , 
 \dot{\gm}_n (u) 
) du 
+ 
\frac{ 
  (m-1) 
  G_n' ( d ( o_n , X^\a ( t_n^{(\a)} ) ) ) 
}{
  G_n ( d ( o_n , X^\a ( t_n^{(\a)} ) ) ) 
}
\\
& \le 
C_1 
\abra{ 
  r_0 + \coth ( C_1 r_0 ) 
}. 
\end{align*}
Hence the conclusion follows from \eq{expectation} 
with $C_0 = C_1 ( 1 + 3 r_0 / 4 + \coth ( C_1 r_0 ) / 2 )$. 
\epf 

In the next step, we will introduce 
a comparison process to give a control 
of the radial process. 
Let us define two functions $\ph$ and $\psi$ 
on $(2r_0 , \infty )$ by 
\begin{align*}
\ph ( r )
& : = 
C_0 + \frac12 \int_0^r b (u) \, du 
,
\\ 
\psi ( r ) 
& :=   
\frac{2}{r - 2r_0}
,
\end{align*}
where $C_0$ is as in \Lem{drift}. 
Let us define a comparison process $\ro^\a (t)$ 
taking values in $[0,\infty)$
inductively by 
\begin{align*} 
\ro^\a (0)
& := 
d_{g(0)} ( o , x_0 ) + 3 r_0 , 
\\
\ro^\a (t) 
& := 
\ro^\a ( t_n^{(\a)} ) 
+ \frac{ t - t_n^{(\a)} }{\a^2} 
\abra{
  \a \lm_{n+1} 
  + \a^2 
  ( 
   \ph ( \ro^\a ( t_n^{(\a)} ) ) 
   +
   \psi ( \ro^\a ( t_n^{(\a)} ) ) 
  )
}, \quad t \in [ t_n^{(\a)} , t_{n+1}^{(\a)} ]. 
\end{align*}
The term 
$\psi ( \ro^\a ( t_n^{(\a)} ) )$ is inserted 
to avoid a difficulty coming from 
the absence of the estimate 
in \Lem{variation} 
on a neighborhood of $o$. 
By virtue of this extra term, 
$\ro^\a (t) > 2 r_0$ holds 
for all $t \in [ T_1 , T_2 ]$ 
if $\a$ is sufficiently small. 
Let $\hat{\sg}_R'$ and $\bar{\sg}_R'$ be 
given by 
$
\hat{\sg}_R' 
: =  
\sg_R ( \ro^\a )
$ and 
$
\bar{\sg}_R' 
: = 
\a^{-2} ( \ebra{ \hat{\sg}_R' }_\a - T_1 ) + 1
$. 
The following is a modification of 
an argument in the proof of 
\cite[Theorem~3.5.3]{Hsu} 
into our discrete setting. 
%
\begin{lem} \label{lem:comparison1}
For $\dl >0$, 
there exist a family of events $( E_\dl^{\a} )_\a$ 
with $\lim_{\a \to 0} \P [ E_\dl^{\a} ] = 1$ 
and a constant $K (\dl) > 0$ 
with $\lim_{\dl \to 0} K (\dl) = 0$ 
such that, on $E_\dl^{\a}$, 
\[ 
d_{ g(t) } ( o, X^\a ( t ) ) 
\le 
\ro^\a (t) + K ( \dl ) 
\]
for $t \in [T_1 , \hat{\sg}_R \wg \hat{\sg}_R' \wg T_2 ]$ 
and sufficiently small $\a$ 
relative to $\dl$ and $R^{-1}$. 
\end{lem}
\bpf
It suffices to show the assertion in the case 
$t = t_n^{(\a)}$ for some $n \in \N_0$. 
Indeed, once we have shown it, 
\Cor{away} \ref{approx_error} yields  
\begin{align*}
d_{g(t)} ( o , X^\a (t) ) 
& \le 
\e^{\kp \a^2} 
\abra{ 
  d_{g( \ebra{t}_\a )} ( o , X^\a ( \ebra{t}_\a ) ) 
  + 
  h (\a) 
}
\\
& \le 
\ro_{ \ebra{ t }_\a }^\a 
+ K (\dl) 
+ ( \e^{\kp \a^2} - 1 ) R 
+ \e^{\kp \a^2} h (\a) 
\\
& \le 
\ro_t^\a  
+ K (\dl) 
+ \a  
+ ( \e^{\kp \a^2} - 1 ) R 
+ \e^{\kp \a^2} h (\a) 
\end{align*} 
for $t \in [ T_1 , \hat{\sg}_R \wg T_2]$. 
Here we used the fact $\ph \ge 0$ and $\psi > 0$. 
Thus the conclusion can be easily deduced. 

For simplicity of notations, we denote 
$d_{g( t_n^{(\a)} )} ( o , X^\a ( t_n^{(\a)} ) )$ 
and 
$\ro^\a ( t_n^{(\a)} )$ 
by $d_n$ and $\ro_n$ 
respectively 
in the rest of this proof. 
Let us define a sequence of 
$\sF_n$-stopping times $S_l$ by 
$S_0 : = 0$ and 
\begin{align*}
S_{2l+1} 
& : = 
\inf 
\bbra{ 
  j \ge S_{2l} 
  \; \left| \; 
      X^\a ( t_j^{(\a)} ) \in B_{r_0}^{(t_j^{(\a)})} (o) 
  \right.
} \wg N^{(\a)} , 
\\ 
S_{2l} 
& : = 
\inf 
\bbra{ 
  j \ge S_{2l-1} 
  \; \left| \; 
      X^\a ( t_j^{(\a)} ) \notin B_{3r_0 / 2}^{(t_j^{(\a)})} (o)
  \right.
} \wg N^{(\a)}
.
\end{align*}
Since $\ro_n > 2r_0$, 
it suffices to show the assertion 
in the case 
$
S_{2l} 
\le n 
< S_{2l+1} \wg \bar{\sg}_R \wg \bar{\sg}_R' 
$ 
for some $l \in \N_0$. 
Now \Lem{variation} and \Lem{drift} imply 
\[
d_{j+1} - \ro_{j+1} 
\le 
d_{j} - \ro_{j} 
+ \a^2 ( \ph ( d_{j} ) - \ph ( \ro_{j} ) ) 
+ \a^2 ( \Lm_{j+1} - \bar{\Lm}_{j+1} ) 
+ o ( \a^2 )  
\] 
for 
$
j \in 
[ 
  S_{2l} , 
  S_{2l+1} \wg \sg_R' \wg \bar{\sg}_R' 
)
$. 
Here we used the fact $\psi > 0$. 
Let $f_\a$ be a $C^2$-function on $\R$ 
satisfying 
\begin{enumerate}
\item \label{flat}
$f_\a |_{( - \infty , -\a )} \equiv 0$, 
$f_\a |_{( \a , \infty )} (x) = x$, 
\item \label{convex}
$f_\a$ is convex, 
\item \label{2nd_error}
$\a^2 \sup_{x \in \R} f_\a'' (x) = O (1)$.  
\end{enumerate} 
For example, 
a function $f_\a$ satisfying these 
conditions is constructed by setting 
\[
\tilde{f} (x)
= 
\int_{-\infty}^x \int_{-\infty}^t 
b \exp \abra{ - \frac{a}{1-s^2} } 1_{(-1,1)} (s)
ds dt ,
\]
where $a,b$ is chosen to satisfy 
\begin{align*}
\int_{-\infty}^1 & 
\exp \abra{ - \frac{ a }{ 1-s^2 } } 1_{(-1,1)} (s)
ds = 1 , 
& 
b \int_{-\infty}^{1} & 
\int_{-\infty}^t 
\exp \abra{ - \frac{ a }{ 1-s^2 } } 1_{(-1,1)} (s)
ds dt 
= 1 
\end{align*} 
and $f_\a ( x ) : = \a \tilde{f} ( \a^{-1} x )$. 
By the Taylor expansion 
with the condition \ref{2nd_error} of $f_\a$, 
we have 
\begin{multline} \label{eq:Taylor}
f_\a ( d_{j+1} - \ro_{j+1} )
 \le 
f_\a ( d_{j} - \ro_{j} ) 
\\
+ 
\a^2 
f_\a' ( d_{j} - \ro_{j} ) 
( 
 \ph ( d_{j} ) - \ph ( \ro_{j} ) 
 + 
 ( \Lm_{j} - \bar{\Lm}_{j} )
) 
+ o ( \a^2 ) . 
\end{multline} 
Let $C > 0$ be the Lipschitz constant of $\ph$ 
on $[ 0, R ]$. 
Note that we have 
\begin{equation} \label{eq:Taylor_1st}
f_\a' ( d_{j} - \ro_{j} ) 
( \ph ( d_{j} ) - \ph ( \ro_{j} ) )
\le C ( d_{j} - \ro_{j} )_+ + o (1) 
. 
\end{equation}
Here the error term $o(1)$ may appear 
in the case $d_{j} - \ro_{j} \in [ -\a , 0]$. 
Now by using \eq{Taylor} and \eq{Taylor_1st} 
combined with the fact 
$d_{S_{2l}} - \ro_{S_{2l}} < - \a$ 
for sufficiently small $\a$, 
we obtain 
\begin{align} \nn
( d_{n} - \ro_{n} )_+ 
& \le 
f_\a ( d_{n} - \ro_{n} ) 
\\ 
\label{eq:pre-comparison} 
& \le 
C \a^2 
\sum_{j = S_{2k}}^{n-1} 
( d_{j} - \ro_{j} )_+ 
+ \a^2 \sum_{j = S_{2k}}^{n-1} 
f_\a' ( d_{j} - \ro_{j} ) 
( \Lm_{j+1} - \bar{\Lm}_{j+1} )
+ o (1) . 
\end{align}
Here the first inequality follows from 
the condition \ref{convex} of $f_\a$ 
and 
$n \le \a^{-2} ( T_2 - T_1 )$ is used 
to derive the error term $o(1)$. 
Let $E_\dl^{\a}$ be an event defined by 
\[
E_\dl^{\a} 
: = 
\bbra{ 
  \a^2 
  \sup_{
    k
     \le 
    k' 
     \le 
    N^{(\a)} \wg \bar{\sg}_R 
  }
  \abs{ 
    \sum_{j=k}^{k'} 
    f_\a' ( d_{j-1} - \ro_{j-1} ) 
    ( \Lm_{j} - \bar{\Lm}_{j} ) 
  } 
  < \dl 
}. 
\]
Note that  
$a_j = f_\a' ( d_{j-1} - \ro_{j-1} )$ is 
$\sF_n$-predictable 
and uniformly bounded by $1$. 
Thus, 
by combining \Lem{LLN} with \eq{pre-comparison}, 
we obtain 
\begin{equation*}
( d_{n} - \ro_{n} )_+ 
\le 
C \a^2 
\sum_{j = S_{2l}}^{n-1} ( d_{j} - \ro_{j} )_+ 
+ 2 \dl 
\end{equation*} 
on $E_\dl^{\a}$ for sufficiently small $\a$. 
Thus, 
by virtue of a discrete Gronwall inequality 
(see \cite{Will-Wong} for instance), 
\[
( d_{n} - \ro_{n} )_+ 
\le 
2 \dl 
\abra{ 
  1 + ( 1 + C \a^2 )^n 
}
\le 
2 \dl ( 1 + \e^{ C (T_2 - T_1)} ) 
. 
\]
This estimate implies the conclusion. 
\epf
\begin{cor} \label{cor:comparison_stopping} 
For every $R' < R$, 
$
\limsup_{\a \to 0 } \P [ \hat{\sg}_R \le T_2 ] 
\le 
\limsup_{\a \to 0} \P [ \hat{\sg}_{R'}' \le T_2 ]
$. 
\end{cor}
Now we turn to the proof of our destination 
in this section. 
%
%
\apf{\Prop{non-explosion}}
By \Cor{comparison_stopping}, 
the proof of \Prop{non-explosion} is reduced 
to estimate 
$\P [ \hat{\sg}_R' \le T_2 ]$. 
To obtain a useful bound of it, 
we would like to apply 
the invariance principle for $\ro^\a$. 
However, there is a technical difficulty 
coming from 
the unboundedness of the drift term of $\ro^\a$. 
To avoid it, 
we introduce an auxiliary process 
$\tilde{\ro}^\a$ in the sequel. 

Let $\tilde{\ph}$ be a bounded, 
globally Lipschitz 
function on $\R$ such that 
$\tilde{ \ph } (r) = \ph (r) + \psi (r)$ 
for $r \in [ 2 r_0 + R^{-1} , R ]$. 
Let us define an $\R$-valued process 
$\tilde{\ro}^\a (t)$ inductively by 
\begin{align*} 
\tilde{\ro}^\a (0) 
& := 
d_{g(T_1)} ( o , x_0 ) + 3 r_0 , 
\\
\tilde{\ro}^\a (t) 
& := 
\tilde{\ro}^\a ( t_n^{(\a)} )
+ \frac{ t - t_n^{(\a)} }{\a^2} 
\abra{
  \a \lm_{n+1}  
  + \a^2 
  \tilde{\ph} ( \ro^\a ( t_n^{(\a)} ) ) 
}, 
\quad 
t \in [ t_n^{(\a)} , t_{n+1}^{(\a)} ] . 
\end{align*}
We also define two diffusion processes 
$\ro^0 (t)$ and $\tilde{\ro}^0 (r)$ 
as solutions to the 
following SDEs:   
\begin{align*}
& \left\{ 
\begin{array}{ll}
d \ro^0 (t) 
& \hspace{-.5em} = 
d B(t) 
+ ( \ph ( \ro^0 (t) ) + \psi ( \ro^0 (t) ) ) dt ,
\\
\ro^0 (T_1) 
& \hspace{-.5em} = 
d_{g(T_1)} ( o , x_0 ) + 3 r_0 , 
\end{array}
\right.
\\
& \left\{ 
\begin{array}{ll}
d \tilde{\ro}^0 (t) 
& \hspace{-.5em} = 
d B(t) 
+ \tilde{\ph} ( \tilde{\ro}^0(t) ) d t , 
\\
\tilde{\ro}^0 (T_1) 
& \hspace{-.5em} = 
d_{g(T_1)} ( o , x_0 ) + 3 r_0 , 
\end{array}
\right. 
\end{align*} 
where $( B(t) )_{t \in [ T_1 , T_2 ]}$ is 
a standard 1-dimensional Brownian motion 
with $B(T_1) = 0$. 
We claim that $\tilde{\ro}^\a$ converges 
in law to $\tilde{\ro}^0$ as $\a \to 0$. 
Indeed, we can easily show 
the tightness of $( \ro^\a )_{\a > 0}$ 
by modifying an argument 
for the invariance principle for 
i.i.d.~sequences since $\tilde{\ph}$ is bounded. 
Then the claim follows 
from the same argument 
as we used in the proof of \Thm{IP} 
under \Prop{tight}, which is 
based on the Poisson subordination and 
the uniqueness of the martingale problem. 

Let us define 
$
\h_R 
\: : \: 
\sC_1 \to [ T_1 , T_2 ] \cup \{ \infty \}
$ 
by 
$
\h_R (w) 
= 
\inf 
\bbra{ 
  t \in [ T_1 , T_2 ]  
  \; | \; 
  w (t) \le 2r_0 + R^{-1} 
} 
$. 
Then we have 
\begin{equation*}
\P [ \hat{\sg}_R' \le T_2 ]
\le 
\P 
\cbra{ 
  \sg_R ( \ro^\a ) \wg \h_R ( \ro^\a ) \le T_2 
} 
=  
\P 
\cbra{ 
  \sg_R ( \tilde{\ro}^\a ) \wg \h_R ( \tilde{\ro}^\a )  
  \le T_2 
}. 
\end{equation*}
Since $\{ w \; | \; \sg_R (w) \wg \h_R (w) \le T_2 \}$ 
is closed in $\sC_1$, 
the Portmanteau theorem implies 
\[
\limsup_{\a \to 0} 
\P 
\cbra{ 
  \sg_R ( \tilde{\ro}^\a ) \wg \h_R ( \tilde{\ro}^\a )  
  \le T_2 
} 
\le 
\P 
\cbra{ 
  \sg_R ( \tilde{\ro}^0 ) \wg \h_R ( \tilde{\ro}^0 )  
  \le T_2 
} 
= 
\P 
\cbra{ 
  \sg_R ( \ro^0 ) \wg \h_R ( \ro^0 )  
  \le T_2 
}. 
\]
Since $\ro^0$ is a diffusion process on $( 2 r_0 , \infty )$ 
which cannot reach the boundary by \Ass{non-explosion}, 
the conclusion follows.  
\epf

\subsection{Tightness of geodesic random walks}
\label{sec:tight}

Recall that we have metrized the path space $\sC$ 
by using $d_{g(T_1)}$. 
To deal with 
the tightness of $( X^\a )_{\a \in ( 0, 1 )}$ 
in $\sC$, 
we show the following lemma, 
which provides a tightness criterion compatible  
with the time-dependent metric $d_{g(t)}$. 
\begin{lem} \label{lem:suff_tight} 
$( X^\a )_{\a \in ( 0, 1 )}$ is tight 
if 
\begin{equation} \label{eq:tight}
\lim_{\dl \to 0} \frac{1}{\dl} 
\limsup_{\a \to \infty} 
\sup_{
      n \in \N_0 
}
\P \cbra{ 
  \sup_{t_n^{(\a)} \le s \le ( t_n^{(\a)}  + \dl ) \wg T_2} 
  d_{g(s)} ( X^\a ( t_n^{(\a)} ) , X^\a (s) ) > \ep , 
  \hat{\sg}_R = \infty 
}
= 0 
\end{equation} 
holds for every $\ep > 0$ and $R > 1$. 
\end{lem} 
\bpf
By following a standard argument 
(cf. \cite[Theorem~7.3 and Theorem~7.4]{Billi}),  
we can easily show that 
$( X^\a )_{\a \in (0,1)}$ is tight 
if, for every $\ep > 0$, 
\begin{equation} \nn 
\lim_{\dl \to 0} \frac{1}{\dl} 
\limsup_{\a \to 0} 
\sup_{t \in [ T_1 , T_2 ]} 
\P 
\cbra{ 
  \sup_{t \le s \le ( t + \dl ) \wg T_2} 
  d_{g(T_1)} ( X^\a (t) , X^\a (s) ) > \ep 
}
= 0
.
\end{equation} 
Thus, by virtue of \Prop{non-explosion}, 
$ ( X^\a )_{\a \in ( 0 , 1 )}$ is tight if 
\begin{equation} \nn 
\lim_{\dl \to 0} \frac{1}{\dl} 
\limsup_{\a \to 0} 
\sup_{t \in [ T_1 , T_2 ]} 
\P 
\cbra{ 
  \sup_{t \le s \le ( t + \dl ) \wg T_2} 
  d_{g(T_1)} ( X^\a (t) , X^\a (s) ) > \ep , 
  \hat{\sg}_R = \infty
}
= 0 
\end{equation} 
for every $\ep > 0$ and $R > 1$.
Given $R > 1$, 
take $M_0$ and $\kp$ 
as in \Lem{bdd_cpt} and \Lem{metric_control}
respectively. 
Then, for $\ep < 1$ and $s,t \in [ T_1 , T_2 ]$, 
\[
\bbra{ 
  d_{g(s)} ( X^\a (s) , X^\a ( \ebra{ t }_\a ) ) 
  \le \ep , 
  \hat{\sg}_R = \infty 
}  
\subset 
\bbra{ 
  d_{g(T_1)} ( X^\a (s), X^\a ( t ) ) 
  \le 2 \e^{\kp ( T_2 - T_1 )} \ep , 
  \hat{\sg}_R = \infty
} 
\] 
if $\a$ is sufficiently small. 
Thus we have 
\begin{multline*}
\bbra{ 
  \sup_{t \le s \le ( t + \dl ) \wg T_2} 
  d_{g(T_1)} ( X^\a (t) , X^\a (s) ) > \ep , 
  \hat{\sg}_R = \infty 
}
\\ 
\subset 
\bbra{ 
  \sup_{\ebra{t}_\a \le s \le ( \ebra{t}_\a + 2 \dl ) \wg T_2} 
  d_{g(s)} ( X^\a ( \ebra{t}_\a ) , X^\a (s) ) 
  > \frac{\e^{-\kp ( T_2 - T_1 )} \ep}{2} , 
  \hat{\sg}_R = \infty 
}
\end{multline*}
for $\a^2 \le \dl$ and hence the conclusion follows. 
\epf
\apf{\Prop{tight}}
Take $R > 1$. 
By virtue of \Lem{suff_tight}, 
it suffices to show \eq{tight}. 
Take $M_0 \subset M$ compact 
and $\kp$ as in \Lem{bdd_cpt} and 
\Lem{metric_control} respectively. 
By taking smaller $\ep > 0$, 
we may assume that $\ep < \tilde{r_0} / 2$, 
where $\tilde{r}_0 = \tilde{r}_0 (M_0)$ 
is as in \Lem{inject}. 
Take $n \in \N_0$ with $n < N^{(\a)}$. 
Let us define a $\sF_k$-stopping time $\zt_\ep$ 
by 
\[
\zt_\ep 
: = 
\inf \bbra{
  k \in \N_0 
  \; | \; 
  n \le k \le N^{(\a)} , 
  d_{g ( t_k^{(\a)} )} ( X^\a ( t_n^{(\a)} ) , X^\a ( t_k^{(\a)} ) ) 
  > \ep 
}. 
\] 
Then, for sufficiently small $\a$, 
\begin{multline} \label{eq:replace}
\bbra{ 
  \sup_{ t_n^{(\a)} \le s \le ( t_n^{(\a)} + \dl ) \wg T_2 } 
  d_{g(s)} ( X^\a ( t_n^{(\a)} ) , X^\a ( s ) ) \ge 2 \ep 
  , \hat{\sg}_R = \infty
}
\\
\subset 
\bbra{ 
  \a^2 ( \zt_{\ep} - n ) < \dl 
  , \, 
  \hat{\sg}_R = \infty 
} . 
\end{multline}
Set $p_k := X^\a ( t_k^{(\a)} )$ for $k \in \N_0$ 
and 
$f ( t , x ) := d_{g(t)} ( p_n , x )$. 
Note that $f^2$ is smooth 
on $\bbra{ f < \ep }$. 
Let us define $\lm_k'$ by 
\[
\lm_{k+1}' 
:= 
\dbra{ 
  \tilde{\xi}_{k+1} , 
  \dot{\gm}_{ p_n p_k }^{(t_k^{(\a)})} 
}_{g (t_k^{(\a)})} . 
\] 
We claim that there exists a constant $C > 0$ 
such that 
\begin{equation} \label{eq:2nd_tight}
f ( t_{k+1}^{(\a)} , p_{k+1} )^2 
\le 
f ( t_k^{(\a)} , p_k )^2 
+ 
2 \a f( t_k^{(\a)} , p_k ) \lm_{k+1}' 
+ 
C \a^2  
\end{equation}
for $k \le \zt_\ep \wg N^{(\a)}$ 
on $\{ \hat{\sg}_R = \infty \}$. 
Indeed, 
in the same way as we did to obtain \eq{2nd}, 
\begin{align} \nn
f ( t_{k+1}^{(\a)} , p_{k+1} )^2 
& \le 
f ( t_k^{(\a)} , p_k )^2 
 + 
2 \a f ( t_k^{(\a)} , p_k ) \lm_{k+1}' 
 + 
\a^2 ( \lm_{k+1}' )^2 
\\ \nn
& \quad 
 + 
2 \a^2 f( t_k^{(\a)} , p_k ) 
\abra{ 
  \partial_t f ( t_k^{(\a)} , p_k ) 
  + 
  \dbra{ 
    Z ( t_k^{(\a)} ) 
    , 
    \dot{\gm}_{p_n p_k}^{(t_n^{(\a)})} 
  }_{g ( t_k^{(\a)} )} ( p_k ) 
} 
\\ \label{eq:2nd_tight1} 
& \quad 
 + 
\a^2 
f( t_k^{(\a)} , p_k ) 
I_{\gm_{p_n p_k}^{(t_k^{(\a)})}}^{( t_k^{(\a)} )} 
( 
 J_{\tilde{\xi}_{k+1}} 
) 
 + 
o ( \a^2 ).  
\end{align}
Here $o ( \a^2 )$ is controlled uniformly. 
Let $K_1 > 0$ be a constant satisfying that 
the $g(t)$-sectional curvature on $M_0$ 
is bounded below by $-K_1$ for every $t \in [ T_1 , T_2 ]$. 
Such a constant exists since $M_0$ is compact. 
Then a comparison argument implies 
\begin{equation*} 
f( t_k^{(\a)} , p_k ) 
I_{\gm_{p_n p_k}^{(t_k^{(\a)})}}^{( t_k^{(\a)} )} 
( 
 J_{\tilde{\xi}_{k+1}} 
) 
\le 
K_1 f ( t_k^{(\a)} , p_k ) 
\coth ( K_1 f ( t_k^{(\a)} , p_k ) ) . 
\end{equation*} 
Here the right hand side is bounded uniformly 
if $k < \zt_\ep \wg N^{(\a)}$. 
The remaining estimate of the second order term in 
\eq{2nd_tight1} to show \eq{2nd_tight} is easy 
since we are on the event $\{ \hat{\sg}_R = \infty \}$. 
Applying \eq{2nd_tight} repeatedly 
from $k=n$ to $k = \zt_\ep$, we obtain 
\begin{align*}
\ep^2 
< 
\a \sum_{k=n}^{\zt_\ep} f ( t_k^{(\a)} , p_k ) \lm_{k+1}' 
+ C \dl 
\end{align*}
on 
$
\bbra{ 
  \a^2 ( \zt_\ep - n ) < \dl 
  , \, 
  \hat{\sg}_R = \infty 
}
$. 
Set
$ 
N_\dl^{(\a)} 
:= 
\sup \{ 
  k \in \N_0 
  \; | \; 
  k \le \a^{-2} \dl + n 
\}
$. 
By taking $\dl < ( 2 C )^{-1} \ep^2$, 
we obtain 
\begin{align} \nn 
\bbra{ 
  \a^2 ( \zt_\ep - n ) < \dl 
  , \, 
  \hat{\sg}_R = \infty 
} 
& \subset 
\bbra{ 
  \sum_{k=n}^{\zt_\ep} 
  f ( t_k^{(\a)} , p_k ) 
  \lm_{k+1}' > \frac{\ep^2}{2 \a}
  , \,
  \a^2 ( \zt_\ep - n ) < \dl 
  , \, 
  \hat{\sg}_R = \infty 
}
\\ \label{eq:replace2}  
& \subset 
\bbra{ 
  \sup_{n \le N \le N_\dl^{(\a)}} 
  \sum_{k=n}^{N} 
  f ( t_k^{(\a)} , p_k ) 
  1_{ \{ f ( t_k^{(\a)} , p_k ) \le \ep \} }
  \lm_{k+1}' > \frac{\ep^2}{2 \a}
}. 
\end{align} 
Set 
\begin{align*}
Y_{k+1} 
& : = 
\frac{1}{\sqrt{m+2}} 
f ( t_k^{(\a)} , p_k ) 
1_{ \{ f ( t_k^{(\a)} , p_k ) \le \ep \} }
\lm_{k+1}' .
\end{align*}
We can easily see that 
$| Y_k | \le 1$ and 
$\sum_{k=n+1}^N Y_k$ is $\sF_N$-martingale. 
By \cite[Theorem~1.6]{Freed_tail} 
with \eq{replace2}, 
we obtain
\begin{align*}
\P [ 
  \a^2 ( \zt_\ep - n ) < \dl
  , \, 
  \hat{\sg}_R = \infty 
] 
& \le 
\P \cbra{ 
\sup_{n \le N \le N_\dl^{(\a)} } 
\sum_{k=n+1}^{N+1} Y_k  
> \frac{\ep^2}{2 \a \sqrt{m+2}} 
}
\\
& \le 
\exp 
\abra{ 
  - 
  \frac
    {\ep^4}
    {
      4 \sqrt{m+2} 
      ( \a \ep^2 + 2 \a^2 \sqrt{m+2} ( N_\dl^{(\a)} - n ) )
    }  
}
\\
& \le 
\exp 
\abra{ 
  - 
  \frac
    {\ep^4}
    {
      4 \sqrt{m+2} 
      ( \a \ep^2 + 2 \sqrt{m+2} \dl )
    }  
}
.
\end{align*}
Hence \eq{tight} follows 
by combining this estimate with \eq{replace}. 
\epf

\section{Coupling by reflection}
\label{sec:reflection}

For $k \in \R$, 
let $U_{a,k}$ be a 1-dimensional 
Ornstein-Uhlenbeck process 
defined as a solution to the following SDE: 
\begin{align*}
d U_{a,k} (t) & = - \frac{k}{2} U_{a,k} (t) dt + 2 d B (t) , 
\\
U_{a,k} ( T_1 ) & = a. 
\end{align*}
More explicitly, 
$
U_{a,k} (t) 
= 
\e^{-k(t - T_1)/2} a  + 2 \int_0^t \e^{k (s-t)/2} d B(s)
$. 
Here $B(t)$ is standard 1-dimensional Brownian motion 
as in the proof of \Prop{non-explosion}. 
\begin{thm} \label{th:main0} 
Suppose 
\begin{equation} \label{eq:c-b0} 
2 ( \nab Z (t) )^\flat + \partial_t g (t) 
\le 
\Ric_{g(t)} + k g(t) 
\end{equation} 
holds for some $k \in \R$. 
Then, for each $x_1 , x_2 \in M$, 
there exists a coupling 
$\mathbf{X} (t) := ( X_1 (t) , X_2 (t) )$ 
of two $\sL_t$-diffusion particles starting at $(x_1 , x_2 )$ 
satisfying 
\begin{equation} \nn 
\P \cbra{ 
  \inf_{T_1 \le t \le T} d_{g(t)} ( \mathbf{X}(t) ) > 0 
} 
\le 
\P \cbra{ 
  \inf_{T_1 \le t \le T}  U_{d_{g(T_1)} (x_1, x_2),k} (t) > 0  
} 
= 
\chi 
\abra{ 
  \frac
  { d_{g(T_1)} (x_1 , x_2) }
  { 2 \sqrt{ \b ( T - T_1 ) } } 
} 
\end{equation} 
for each $T \in [ T_1 , T_2 ]$, 
where 
\begin{align*}
\chi (a) 
& := 
\frac{1}{\sqrt{2 \pi}} 
\int_{-a}^a \e^{-u^2 /2} du, 
& 
\b (t) 
& := 
\begin{cases} 
\ds 
\frac{e^{kt} - 1}{k} 
& k \neq 0 ,
\\ \ds
t 
& k = 0 . 
\end{cases}
\end{align*}
\end{thm}
\begin{remark} \label{rem:O-U}
\begin{enumerate}
\item
Given $k \in \R$, 
a simple example satisfying \eq{c-b0} can be constructed 
from a solution $\tilde{g} (t)$ to the Ricci flow 
$\partial_t \tilde{g} (t) = \Ric_{\tilde{g}(t)}$ by a scaling. 
That is, $g(t) = \e^{-kt} \tilde{g} (t)$ satisfies \eq{c-b0} 
(with equality) when $Z(t) \equiv 0$. 
\item
Our assumption \eq{c-b0} can be regarded 
as a natural extension of a lower Ricci curvature bound 
by $k$. 
Indeed, Bakry-\'{E}mery's curvature-dimension condition 
${\sf CD} (k ,\infty)$  (see \cite{Bak97} for instance), 
which is a natural extension of a lower Ricci curvature bound, 
appears in \eq{c-b0} 
when both $Z(t)$ and $g(t)$ are independent of $t$. 
\item 
From the last item in this remark, 
when $Z(t) \equiv 0$, 
one may expect that \eq{c-b0} works 
as an analogue of Bakry-\'{E}mery's 
${\sf CD} (k , N)$ condition, 
which is equivalent to $\Ric_g \ge k$ and $\dim M < N$ 
when $g(t)$ is independent of $t$, 
instead of ${\sf CD}(k, \infty)$ 
since $\dim M = m < \infty$ in our case. 
However, 
the following observation suggests us that 
we should be more careful: 
Let us consider \eq{c-b0} 
in the case $k > 0$ and $Z(t) \equiv 0$. 
When $\partial_t g(t) \equiv 0$, 
the Bonnet-Myers theorem tells us that 
the diameter of $M$ is bounded and 
hence $M$ is compact. 
When $g(t)$ depends on $t$, it is no longer true. 
In fact, we can easily obtain a noncompact $M$ 
satisfying \eq{c-b0} with $k > 0$ 
by following an observation 
in the first item of this remark. 
On the other hand, 
the Bonnet-Myers theorem is known to hold 
under ${\sf CD} (k, N)$ when $Z$ is of the form $\nab h$ 
in the time-homogeneous case 
(see \cite{Bak-Led_BM-Sob,Qian:BonnetMyers}). 
\end{enumerate}
\end{remark} 
By following a standard argument, 
\Thm{main0} implies the following estimate 
for a gradient of the diffusion semigroup: 
\begin{cor} \label{cor:sF}
Let $( \{ X(t) \}_{t \in [ T_1, T_2 ]} ,  \{ \P_x \}_{x \in M} )$ 
be a $\sL_t$-diffusion process with $\P_x [ X(T_1) = x ] = 1$. 
For any bounded measurable function $f$ on $M$, 
let us define $P_t f$ by $P_t f (x) := \E_x [ f ( X(t) ) ]$. 
Then, under the same assumption as in \Thm{main0}, 
we have 
\[
\limsup_{y \to x} 
\abs{ \frac{ P_t f (x) - P_t f (y) }{d_{g(T_1)} (x,y)} }
\le 
\frac{1}{\sqrt{2 \pi \b ( t - T_1 )}} 
\sup_{z , z' \in M} | f (z) - f (z') | 
.
\]
In particular, 
$P_t f$ is 
$d_{g(T_1)}$-globally Lipschitz continuous 
when $f$ is bounded. 
\end{cor}
\apf{\Cor{sF}}
Let $\mathbf{X} = (X_! , X_2 )$ be 
a coupling of $\sL_t$-diffusions 
$( X(t) , \P_x )$ and $( X(t) , \P_y )$ 
given in \Thm{main0}. 
Let $\tau^*$ be 
the coupling time of $\mathbf{X}$, 
i.e. 
$
\tau^* 
: = 
\inf \{ 
  t \in [ T_1 , T_2 ] 
  \; | \; 
  \mathbf{X} (t) \in D (M) 
\}
$. 
Let us define a new coupling 
$\mathbf{X}^* = ( X_1^* , X_2^* )$ 
by 
\begin{equation} \nn 
\mathbf{X}^* (t) : = 
\begin{cases} 
\mathbf{X} (t) & \mbox{if $\tau^* > t$,} 
\\
( X_1 (t) , X_1 (t) ) & \mbox{otherwise.}
\end{cases}
\end{equation}
Since 
$
\{ \tau^* > T \} 
= 
\{ \inf_{T_1 \le t \le T} d_{g(t)} ( \mathbf{X} (t) ) > 0 \}
$, 
\Thm{main0} yields 
\begin{align} \nn
P_t f (x) - P_t f (y) 
& = 
\E [ f ( X_1^* (t) ) - f ( X_2^* (t) ) ] 
\\ \nn
& = 
\E \cbra{ 
  \abra{ 
    f ( X_1^* (t) ) - f ( X_2^* (t) ) 
  }
  1_{ \{ \tau^* > t \} } 
} 
\\ \nn
& \le 
\P [ \tau^* > t ] \sup_{z,z' \in M} | f (z) - f(z') | 
\\ \nn 
& \le 
\chi \abra{ \frac{d_{g(T_1)} (x,y)}{2 \sqrt{\b ( t - T_1 )}} }
\sup_{z,z' \in M} | f(z) - f(z') |
. 
\end{align}
Hence the assertion holds 
by dividing the both sides of the above inequality 
by $d_{g(T_1)} (x,y)$ 
and 
by letting $y \to x$ after that. 
\epf

%
As we did in the last section, 
let $( \gm_{xy}^{(t)} )_{x,y \in M}$ be 
a measurable family of 
unit-speed minimal $g(t)$-geodesics 
such that $\gm_{xy}^{(t)}$ joins $x$ and $y$. 
Without loss of generality, 
we may assume that 
$\gm_{xy}^{(t)}$ is symmetric, that is, 
$
\gm_{xy}^{(t)} ( d_{g(t)} (x,y) - s ) 
= 
\gm_{yx}^{(t)} (s)
$ 
holds. 
Let us define 
$\tilde{m}_{xy}^{(t)} \: : \: T_y M \to T_y M$ 
by 
\[
\tilde{m}_{xy}^{(t)} v 
:= 
v 
- 2 \dbra{ v, \dot{\gm}_{xy}^{(t)} }_{g(t)} 
\dot{\gm}_{xy}^{(t)} ( d_{g(t)} (x,y) )
. 
\] 
This is a reflection with respect to a hyperplane 
which is $g(t)$-perpendicular to $\dot{\gm}_{xy}^{(t)}$. 
Let us define $m_{xy}^{(t)} \: : \: T_x M \to T_y M$ by 
$m_{xy}^{(t)} := \tilde{m}_{xy}^{(t)} \circ /\!\!/_{\gm_{xy}^{(t)}}^{(t)}$. 
Clearly $m_{xy}^{(t)}$ is a $g(t)$-isometry. 
As in the last section, 
let $\Ph^{(t)} \: : \: M \to \mathscr{O}^{(t)} (M)$ be 
a measurable section of 
the $g(t)$-orthonormal frame bundle 
$\mathscr{O}^{(t)} (M)$ of $M$. 
Let us define two measurable maps 
$\Phi_i^{(t)} \: : \: M \times M \to \mathscr{O}^{(t)} (M)$ 
for $i=1,2$ by 
\begin{align*}
\Ph_1^{(t)} (x,y) & 
:= \Ph^{(t)} (x), 
\\ 
\Ph_2^{(t)} (x,y) 
& := 
\begin{cases}
m_{xy}^{(t)} \Ph_1^{(t)} (x,y), 
& 
(x,y) \in M \times M \setminus D(M), 
\\ 
\Phi^{(t)} (x), 
& 
(x,y) \in D (M).
\end{cases}
\end{align*}
Take $x_1 , x_2 \in M$. 
By using $\Phi_i^{(t)}$, 
we define a coupled geodesic random walk 
$\mathbf{X}^\a (t) = ( X_1^\a (t) , X_2^\a (t) )$ 
by $X^\a_i (0) = x_i$ and, 
for $t \in [ t_n^{(\a)} , t_{n+1}^{(\a)} ]$, 
\begin{align} 
\nn
\tilde{\xi}_{n+1}^i 
& : = 
\sqrt{m+2} 
\Phi_i^{( t_n^{(\a)} )} 
\abra{ 
  \mathbf{X}^\a ( t_n^{(\a)} ) 
} 
\xi_{ n + 1 } ,
\\ \nn 
X_i^\a ( t ) 
& := 
\exp_{X_i^\a ( t_n^{(\a)} )}^{( t_n^{(\a)} )} 
\bigg(
  \frac{ t - t_n^{(\a)} }{\a^2}  
  \Big( 
    \a \tilde{\xi}_{n+1}^i 
     + 
    \a^2 Z( t_n^{(\a)} )
  \Big)
\bigg)  
\end{align}
for $i = 1, 2$. 
We can easily verify that $X_i^\a$ has 
the same law as $X^\a$ with $x_0 = x_i$. 
%
%

In what follows, we assume \eq{c-b0}. 
We can easily verify that 
it implies \Ass{non-explosion}. 
Thus, by \Thm{IP}, 
$( \mathbf{X}^\a )_{\a > 0}$ is tight 
under \Ass{non-explosion}. 
In addition, a subsequential limit 
$\mathbf{X}^{\a_k} \to \mathbf{X} = ( X_1 , X_2 )$ 
in law exists and it is a coupling of 
two $\sL_t$-diffusion processes 
starting at $x_1$ and $x_2$ respectively. 
We fix such a subsequence $( \a_k )_{k \in \N}$. 
In the rest of this paper, 
we use the same symbol $\mathbf{X}^\a$ 
for the subsequence $\mathbf{X}^{\a_k}$ 
and the term ``$\a \to 0$'' always means 
the subsequential limit ``$\a_k \to 0$''. 
Set 
$
 \hat{\sg}_R^{i} 
 : = 
 \sg_R ( d_{g(\cdot)} ( o , X_i^\a (\cdot) ) )
$
for $i = 1,2$. 
We fix $R > 1$ sufficiently large 
until the beginning of the proof of \Thm{main0}. 
Let $M_0 \subset M$ be a relatively compact open set 
satisfying \eq{bdd_cpt} for $2R$ instead of $R$. 

We first show a difference inequality of 
$d_{g(t)} ( \mathbf{X}^\a (t) )$. 
To describe it, 
we will introduce several notations 
as in the last section. 
For simplicity, let us denote 
$\gm_{X_1^\a ( t_n^{(\a)} ) X_2^\a ( t_n^{(\a)} )}^{( t_n^{(\a)} )}$ 
by $\bar{\gm}_n$. 
Let us define a vector field $V_{n+1}$ along 
$\bar{\gm}_n$ by 
\[
V_{n+1} 
: = 
/\!\!/_{\bar{\gm}_n}^{(t_n^{(\a)})} 
\abra{  
  \tilde{\xi}_{n+1}^1 
  - 
  \dbra{ 
    \tilde{\xi}_{n+1}^1 , 
    \dot{\bar{\gm}}_n  
  }_{g (t_n^{(\a)})} \dot{\bar{\gm}}_n (0)
}
.
\]
Take $v \in \R^m$. 
Let us define $\lm_{n+1}^*$ and $\Lm_{n+1}^*$ by 
\begin{align*}
\lm_{n+1}^* 
& : = 
\begin{cases}
2 
\dbra{ 
  \tilde{\xi}_{n+1}^1 , 
  \dot{\bar{\gm}}_n 
}_{g(t_n^{(\a)})} 
& \mbox{if $(y_1 , y_2) \notin D(M)$},
\\
2 \sqrt{m+2} 
\dbra{ \xi_{n+1} , v }
& \mbox{otherwise},
\end{cases}
\\
\Lm_{n+1}^* 
& : = 
\frac12 \Bigg( 
\int_0^{d_{g(t_n^{(\a)})} ( \mathbf{X}^\a (t_n^{(\a)}) )} 
\abra{ 
  \partial_t g ( t_n^{(\a)} ) 
  + 
  2 ( \nab Z ( t_n^{(\a)} ) )^\flat 
}
\abra{ 
  \dot{\bar{\gm}}_n (s) , 
  \dot{\bar{\gm}}_n (s) 
}
ds   
\\ 
& \hspace{18em} + 
I_{\bar{\gm}_n}^{(t_n^{(\a)})} 
\abra{ V_{n+1} } 
\Bigg) 1_{\{ \mathbf{X}^\a ( t_n^{(\a)} ) \notin D(M) \} }. 
\end{align*}
For $\dl \ge 0$, let us define  
$\tau_\dl \: : \: \sC_1 \to [ T_1 , T_2 ] \cup \{ \infty \}$ 
by $\tau_\dl (w) := \inf \bbra{ t \ge T_1 \; | \; w(t) \le \dl }$. 
We also define  
$
\hat{\tau}_\dl 
$
by 
$
\hat{\tau}_\dl 
: = 
\tau_\dl ( d_{g(\cdot)} ( \mathbf{X}^a (\cdot) ) )
$. 
%
\begin{lem} \label{lem:c-2var}
For $n \in \N_0$ with $n < N^{(\a)}$, 
we have 
\begin{align} \nn
\e^{ k t_{n+1}^{(\a)} / 2} 
d_{g( t_{n+1}^{(\a)} )} ( \mathbf{X}^\a ( t_{n+1}^{(\a)} ) ) 
& \le 
\abra{ 1 + \frac{k}{2} }
\e^{ k t_n^{(\a)} / 2} 
d_{g( t_n^{(\a)} )} ( \mathbf{X}^\a ( t_n^{(\a)} ) ) 
\\ \label{eq:c-2var}
& \qquad + 
\e^{ k t_n^{(\a)} / 2} ( \a \lm_{n+1}^* + \a^2 \Lm_{n+1}^* ) + o (\a^2 ) 
\end{align}
when 
$n < \hat{\tau}_\dl \wg \hat{\sg}_R^1 \wg \hat{\sg}_R^2$ 
and $\a$ is sufficiently small. 
Moreover, we can control the error term $o(\a^2)$ 
uniformly in the position of $\mathbf{X}^\a$. 
\end{lem}
\bpf
When $( t_n^{(\a)} , \mathbf{X}^\a ( t_n^{(\a)} ) ) \notin \Cutst$, 
\eq{c-2var} is just a consequence of 
the second variational formula for the distance function 
combined with the index lemma for $I_{\bar{\gm}_n}^{(t_n^{(\a)})}$. 
To include the case  
$( t_n^{(\a)} , \mathbf{X}^\a ( t_n^{(\a)} ) ) \in \Cutst$ and 
to obtain a uniform control of $o (\a)$, 
we extend this argument. 
Let us define $H$ and 
$p_1 , p_2 \: : \: H \to [ T_1 , T_2 ] \times \bar{M}_0 \times \bar{M}_0$ 
by 
\begin{align*}
H 
& : = 
\bbra{ 
  (t,x,y,z) 
  \; \left| \; 
  \begin{array}{l}
    t \in [ T_1 , T_2 ], 
    x,y,z \in \bar{M}_0 ,
    \\
    d_{g(t)} (x,y) \ge \dl , 
    \\
    d_{g(t)} (x,y) = 2 d_{g(t)} (x,z) = 2 d_{g(t)} (y,z) 
  \end{array}
  \right.
},
\\
p_1 & (t,x,y,z) 
: = 
(t,x,z), 
\\
p_2 & (t,x,y,z) 
: = 
(t,y,z). 
\end{align*} 
If $\mathbf{q} = (t,x,y,z) \in H$, 
then 
$p_1 (\mathbf{q} ) , p_2 ( \mathbf{q} ) \notin \Cutst$ 
since $z$ is on a midpoint of 
a minimal $g(t)$-geodesic joining $x,y$.    
Since $H$ is compact, 
$p_1 (H)$ and $p_2 (H)$ are also compact. 
Hence there is a constant $\h > 0$ such that 
\[
\inf 
\bbra{ 
  | t - t' | + d_{g(t)} ( x, x' ) + d_{g(t)} ( y, y' ) 
  \; \left| \; 
  \begin{array}{l}
    (t,x,y) \in p_1 (H) \cup p_2 (H) , 
    \\
    (t',x',y') \in \Cutst 
  \end{array}
  \right.
}
> \h . 
\] 
Take $\a > 0$ sufficiently small 
relative to $\h$ and $\dl$. 
Set 
\begin{align*}
p_n 
& := 
\bar{\gm}_n \abra{ 
  \frac{d_{g ( t_n^{(\a)} )} ( \mathbf{X}^\a ( t_n^{(\a)} ) ) }{2}  
}, 
\\
p_n '
& := 
\exp_{p_n}^{t_n^{(\a)}} 
\abra{ 
  V_{n+1} \abra{ 
    \frac{d_{g( t_n^{(\a)} ) } ( \mathbf{X}^\a ( t_n^{(\a)} ) ) }{2}
  }
}
.
\end{align*}
By the triangle inequality, 
we have 
\begin{align*}
d_{g ( t_{n}^{(\a)} )} ( \mathbf{X}^\a ( t_{n}^{(\a)} ) )
& =  
d_{g( t_{n}^{(\a)} )} 
\abra{ 
  X_1^\a ( t_{n}^{(\a)}  ) , 
  p_n 
} 
+ 
d_{g( t_{n}^{(\a)} )} 
\abra{ 
  p_n , 
  X_2^\a ( t_{n}^{(\a)} ) 
} ,    
\\
d_{g ( t_{n+1}^{(\a)} )} ( \mathbf{X}^\a ( t_{n+1}^{(\a)} ) )
& \le 
d_{g( t_{n+1}^{(\a)} )} 
\abra{ 
  X_1^\a ( t_{n+1}^{(\a)}  ) , 
  p_n' 
} 
+ 
d_{g( t_{n+1}^{(\a)} )} 
\abra{ 
  p_n' , 
  X_2^\a ( t_{n+1}^{(\a)} ) 
}.   
\end{align*}
Since 
$
(
 t_n^{(\a)} , 
 \mathbf{X}^\a ( t_n^{(\a)} ), 
 \bar{\gm}_n ( d_{g( t_n^{(\a)} )} ( \mathbf{X}^\a ( t_n^{(\a)} ) ) / 2 ) 
)
\in H
$, 
we can apply the second variational formula 
to each term on the right hand side of 
the above inequality. 
Hence we obtain \eq{c-2var}. 
For a uniform control of the error term, 
we remark that $\bar{\gm}_n$ is included in $M_0$ 
and the $g(t_n^{(\a)})$-length of $\bar{\gm}_n$ is 
bigger than $\dl$. 
These facts follows from 
$n < \hat{\tau}_\dl \wg \hat{\sg}_R^1 \wg \hat{\sg}_R^2$ 
and the choice of $M_0$. 
Thus the every calculation of the second variational 
formula above is done 
on a compact subset of $[ T_1 , T_2 ] \times M_0 \times M_0$ 
which is uniformly away from $\Cutst$. 
It yields the desired result.
\epf

Let us define 
a continuous stochastic process $U_a^\a$ 
on $\R$ starting at $a$ 
by 
\[ 
U_a^\a (t) 
:= 
\e^{-kt/2} a 
+ 
\a \e^{-kt/2} 
\abra{ 
  \sum_{j=1}^{n} 
  \e^{k t_j^{(\a)} /2} \lm_j^* 
  + 
  \frac{ t - t_n^{(\a)} }{\a^2} 
  \e^{ k t_n^{(\a)} /2} 
  \lm_{n + 1}^* 
}
. 
\]
We next show the following comparison theorem 
for the distance process of coupled geodesic random walks. 
\begin{lem} \label{lem:chain-LLN}
For each $\ep > 0$, there exists a family of events 
$( E_\ep^\a )_{\a}$ such that 
$
\P \cbra{ E_\ep^\a }  
$ 
converges to 1
as $\a \to 0$
and 
\begin{equation} \label{eq:chain}
d_{g(t)} ( \mathbf{X}^\a (t) ) 
\le 
U^\a_{d_{g(T_1)} ( \mathbf{X}^\a (T_1) )} (t) 
 +
\ep 
\end{equation} 
for all 
$
t \in 
[ 
  T_1 , 
  T_2 \wg  
  \hat{\tau}_\dl \wg 
  \hat{\sg}_R^1 \wg 
  \hat{\sg}_R^2 
]
$ 
on $E_\ep^\a$ 
for sufficiently small $\a$. 
\end{lem}
\bpf
In a similar way as in the proof of \Lem{comparison1}, 
we can complete the proof once we have found 
$E_\ep^\a$ on which \eq{chain} holds when  
$
t = t_n^{(\a)} 
\in 
[ 
  T_1 , 
  T_2 \wg  
  \hat{\tau}_\dl \wg 
  \hat{\sg}_R^1 \wg 
  \hat{\sg}_R^2 
]
$. 
Set 
$
\bar{\Lm}_{n+1}^*
 := 
\E [ \Lm_{n+1}^* \: | \: \sF_n ]
$ 
and $\bar{\Lm}_0^* := 0$.
Then $\sum_{j=1}^n \e^{ k t_{j-1}^{(\a)} /2 } ( \Lm_j^* - \bar{\Lm}_j^* )$ is 
an $\sF_n$-local martingale. 
Indeed, $\Lm_{n+1}^*$ is bounded 
if $n < \hat{\sg}_R^1 \wg \hat{\sg}_R^2$ and 
so is $\bar{\Lm}_{n+1}^*$. 
Let us define $E_\ep^\a$ by 
\[
E_\ep^\a 
: =  
\bbra{ 
  \sup_{
    \begin{subarray}{c}
        N \le N^{(\a)}  
        \\
        t_N^{(\a)} \le 
        T_2 \wg 
        \hat{\sg}_R^1 \wg 
        \hat{\sg}_R^2 
    \end{subarray}
  }
  \sum_{j=1}^{N+1} 
  \e^{k t_j^{(\a)} / 2} 
  \abra{ \Lm_j^* - \bar{\Lm}_j^* } 
  \le  
  \frac{\ep}{2 \a^2} 
}
. 
\]
In a similar way as in 
\Lem{LLN} or \cite[Lemma~6]{K8}, 
we can show 
$\lim_{\a \to 0} \P \cbra{ E_\ep^\a } = 1$. 
Since we have  
$
(m+2) \E 
[ 
 \dbra{ \xi_i , e_k } 
 \dbra{ \xi_i , e_l } 
] 
= \dl_{kl}
$, 
we obtain 
\[ 
\bar{\Lm}_{n+1}^* 
\le - \frac{k}{2} d_{g( t_n^{(\a)} )} ( \mathbf{X}^\a ( t_n^{(\a)} )). 
\] 
Thus an iteration of \Lem{c-2var} implies 
\eq{chain} on $E_\ep^\a$
when $t = t_n^{(\a)}$. 
\epf


\apf{\Thm{main0}} 
Take $\ep \in ( 0 , 1 )$ arbitrary. 
Let $R > 1$ be sufficiently large 
so that 
\[
\limsup_{\a \to 0} 
\P \cbra{ 
  \hat{\sg}_R^1 \wg \hat{\sg}_R^2 \le T_2 
} 
< \ep 
. 
\]
It is possible by \Prop{non-explosion}. 
Set $a := d_{g(T_1)} ( x_1, x_2 )$. 
Take $T \in [ T_1 , T_2 ]$ and 
let $\dl > 0$ be $\dl >  2 \ep$. 
Then \Lem{chain-LLN} yields 
\begin{align*}
\P \cbra{ 
  \hat{\tau}_\dl > T 
} 
& \le  
\P \cbra{ 
  \bbra{ 
    \hat{\tau}_\dl > T 
  }
  \cap 
  E_\ep^\a 
  \cap 
  \bbra{ 
    \hat{\sg}^1_R \wg \hat{\sg}_R^2 > T 
  } 
} 
+ 2 \ep 
\\
& \le 
\P \cbra{ 
  \tau_{\dl / 2} ( U^\a_a ) > T 
} 
+ 2 \ep 
. 
\end{align*}
Thus we obtain 
\[
\P \cbra{ \hat{\tau}_\dl > T } 
\le 
\P \cbra{ \inf_{t\in [ T_1 ,T ]} U^\a_a (t) \ge \dl / 2 }
\] 
by letting $\ep \downarrow 0$. 
Note that $U^\a_a$ converges in law 
to $U_a$ as $\a \to 0$. 
Since 
\[
\{ 
  \mathbf{w} \in C ( [ T_1,  T_2 ] \to M \times M )
  \; | \; 
  \tau_\dl (d_{g(\cdot)} ( \mathbf{w} (\cdot) ) ) > T 
\}
\] 
is open 
and 
$\bbra{ w \; | \; \inf_{t\in [ T_1 ,T_2 ]} w(t) \ge \dl / 2 }$ 
is closed in $C ( [ 0, T ] \to \R )$,  
the Portmanteau theorem yields 
\begin{align*}
\P \cbra{ 
  \inf_{T_1 \le t \le T} 
  d_{g(t)}( \mathbf{X}( t ) )  
   > 
  \dl 
} 
& \le 
\liminf_{\a \to 0}
\P \cbra{ \hat{\tau}_\dl > T } 
\\
& \le 
\limsup_{\a \to 0} 
\P \cbra{ \inf_{t\in [ T_1 ,T ]} U^\a_a (t) \ge \dl / 2 }
\le
\P \cbra{ \inf_{t\in [ T_1 ,T ]} U_a (t) \ge \dl / 2 }
. 
\end{align*} 
Therefore the conclusion follows 
by letting $\dl \downarrow 0$.     
\epf 
We can also construct a coupling 
by parallel transport by following our manner. 
In the construction of the coupling by reflection, 
we used a map $m_{xy}^{(t)}$. 
By following the same argument 
after replacing $m_{xy}^{(t)}$ with 
$/\!\!/_{\gm_{xy}^{(t)}}^{(t)}$, 
we obtain a coupling by parallel transport. 
The difference of it 
from the coupling by reflection 
is the absence of the term corresponding to $\lm_n^*$, 
which comes from the first variation of arclength. 
As a result, we can show the following 
(cf. \cite{K8}): 
\begin{thm} \label{th:parallel}
Assume \eq{c-b0}. 
For $x_1 , x_2 \in M$, there is 
a coupling $\mathbf{X}(t) = ( X_1 (t) , X_2 (t) )$ 
of two $\sL_t$-diffusion particles 
starting at $x_1$ and $x_2$ at time $T_1$ respectively 
such that 
\[
d_{g(t)} ( \mathbf{X} (t) ) 
\le 
\e^{-k(t-s)/2} 
d_{g(s)} ( \mathbf{X} (s) )
\]
for $T_1 \le s \le t \le T_2$ almost surely. 
\end{thm}
It recovers a part of results 
studied in \cite{Arn-Coul-Thal_horiz}. 
In particular, 
a contraction type estimate 
for Wasserstein distances 
under the heat flow follows. 
\bpf
Let us construct a coupling 
by parallel transport of geodesic random walks 
$
\mathbf{X}^\a 
= 
( X_1^\a  , X_2^\a )
$ 
starting at $(x_1 , x_2 ) \in M \times M$ 
by following the procedure 
stated just before \Thm{parallel}. 
By taking a subsequence, 
we may assume that 
$\mathbf{X}^\a$ converges in law as $\a \to 0$. 
We denote the limit by  $\mathbf{X} = ( X_1 , X_2 )$. 
In what follows, we prove 
\begin{equation*} 
\P \cbra{ 
  \sup_{T_1 \le s \le t \le T_2} 
  \abra{ 
    \e^{kt/2} d_{g(t)} ( \mathbf{X} (t) ) 
    - 
    \e^{ks/2} d_{g(s)} ( \mathbf{X} (s) ) 
  }
  > \ep 
}
= 0  
\end{equation*}
for any $\ep > 0$. 
By virtue of the Portmanteau theorem 
together with \Prop{non-explosion}, 
it suffices to show 
\begin{equation} \label{eq:para0}
\lim_{\a \to 0} 
\P \cbra{ 
  \sup_{T_1 \le s \le t \le T_2} 
  \abra{ 
    \e^{kt/2} d_{g(t)} ( \mathbf{X}^\a (t) ) 
    - 
    \e^{ks/2} d_{g(s)} ( \mathbf{X}^\a (s) ) 
  }
  > \ep 
  , \, 
  \hat{\sg}_R^1 \wg \hat{\sg}_R^2 = \infty 
}
= 0   
\end{equation}
for any $R > 1$. 
For simplicity of notations, 
we write 
$
d_n 
:= 
\e^{k t_n^{(\a)} / 2} d_{g(t_n^{(\a)})} ( \mathbf{X}^\a (t_n^{(\a)}) ) 
$ in this proof. 
For $\dl > 0$, let us define a sequence of 
$\sF_n$-stopping times $S_l$ by 
$S_0 : = 0$ and 
\begin{align*}
S_{2l+1} 
& : = 
\inf 
\bbra{ 
  j \ge S_{2l} 
  \; \left| \; 
      d_j \le \dl 
  \right.
} \wg N^{(\a)} , 
\\ 
S_{2l} 
& : = 
\inf 
\bbra{ 
  j \ge S_{2l-1} 
  \; \left| \; 
      d_j \ge 2 \dl 
  \right.
} \wg N^{(\a)}
.
\end{align*}
Note that $d_{S_{2l-1}} \le 3\dl$ holds 
on $\{ \hat{\sg}_R^1 \wg \hat{\sg}_R^2 = \infty \}$ 
for sufficiently small $\a$. 
As mentioned 
just before \Thm{parallel}, 
\Lem{c-2var} holds with $\lm^* = 0$. 
Moreover, we can obtain the same estimate \eq{c-2var}
even when 
$
S_{2l-1} 
 \le 
n 
 < 
S_{2l} 
\wg \bar{\sg}_R^1
\wg \bar{\sg}_R^2
$ 
for some $l \in \N_0$. 
In this case, the error term $o( \a^2 )$ is 
controlled uniformly also in $l$. 
Let us define an event $E^\a_\dl$ by 
\[
E_\dl^\a 
: =  
\bbra{ 
  \sup_{
    \begin{subarray}{c}
        n \le N \le N^{(\a)}  
        \\
        t_N^{(\a)} \le 
        T_2 \wg 
        \hat{\sg}_R^1 \wg 
        \hat{\sg}_R^2 
    \end{subarray}
  }
  \sum_{j=n+1}^{N+1} 
  \e^{k t_j^{(\a)} / 2} 
  \abra{ \Lm_j^* - \bar{\Lm}_j^* } 
  \le  
  \frac{\dl}{2 \a^2} 
}
. 
\]
Then, as in \Lem{LLN} and \Lem{chain-LLN}, 
we can show $\lim_{\a \to 0} \P [ E_\dl^\a ] = 1$. 
On $E_\dl^\a \cap \{ \hat{\sg}_R^1 \wg \hat{\sg}_R^2 = \infty \}$, 
we have $d_N \le d_n + \dl$ 
for $S_{2l-1} \le n \le N \le S_{2l}$ 
if $\a$ is sufficiently small. 
Moreover, for $n <  S_{2l-1} \le N < S_{2l}$, 
\begin{equation*}
d_N - d_n 
\le 
( d_N - d_{S_{2l-1}} ) + d_{S_{2l-1}} 
\le 
5 \dl . 
\end{equation*}
In the case $S_{2l} \le N < S_{2l+1}$, 
we obtain $d_N - d_n \le 2 \dl$. 
Thus $d_N - d_n \le 5 \dl$ holds for all $n < N$ 
on $E_\dl^\a \cap \{ \hat{\sg}_R^1 \wg \hat{\sg}_R^2 = \infty \}$.  
Take $\dl > 0$ less than $\ep / 10$. 
Then our observations yield \eq{para0} 
since 
$
d_{g(t)} ( \mathbf{X}^\a (t) ) 
 - 
d_{g( \ebra{t}_\a )} ( \mathbf{X} ( \ebra{t}_\a ) )
$ 
becomes uniformly small 
on $\{ \hat{\sg}_R^1 \wg \hat{\sg}_R^2 = \infty \}$ 
as $\a \to 0$.  
\epf



\begin{thebibliography}{10}

\bibitem{Arn-Coul-Thal_horiz}
M.~Arnaudon, K.A. Coulibaly, and A.~Thalmaier, \emph{Horizontal diffusion in
  {$C^1$}-path space}, To appear in S\'{e}minaire de Probabilit\'{e}s, Lecture
  Notes in Mathematics (2009); arXiv:0904.2762.

\bibitem{Bak97}
D.~Bakry, \emph{On {S}obolev and logarithmic {S}obolev inequalities for markov
  semigroups}, New trends in stochasitic analysis (Charingworth, 1994), World
  Sci. Publ. River Edge, NJ, 1997, pp.~43--75.

\bibitem{Bak-Led_BM-Sob}
D.~Bakry and M.~Ledoux, \emph{Sobolev inequalities and {M}yers's diameter
  theorem for an abstract {M}arkov generator}, Duke Math. J. \textbf{85}
  (1996), 253--270.

\bibitem{Billi}
P.~Billingsley, \emph{Convergence of probability measures}, second ed., A Wiley
  Interscienece Publication, John Wiley \& Sons Inc., New York, 1999.

\bibitem{Blum}
Gilles Blum, \emph{A note on the central limit theorem for geodesic random
  walks}, Bull. Austral. Math. Soc. (1984), 169--173.

\bibitem{Chavel2}
I.~Chavel, \emph{Riemannian geometry: a modern introduction}, Cambridge tracts
  in mathematics, 108, Cambridge university press, Cambridge, 1993.

\bibitem{Coul_gtBM}
K.A. Coulibaly, \emph{Brownian motion with respect to time-changing
  {R}iemannian metrics, applications to {R}icci flow}, preprint;
  arXiv:0901.1999.

\bibitem{Crans}
M.~Cranston, \emph{Gradient estimates on manifolds using coupling}, J. Funct.
  Anal. \textbf{99} (1991), no.~1, 110--124.

\bibitem{Ethier-Kurtz}
T.G. Ethier, S.N.~Kurtz, \emph{Markov processes: {C}haracterization and
  convergence.}, Wiley, New York, 1986.

\bibitem{Freed_tail}
D.A. Freedman, \emph{On tail probabilities for martingales}, Ann. Probab.
  \textbf{3} (1975), 100--118.

\bibitem{Hsu}
E.~P. Hsu, \emph{Stochastic analysis on manifolds}, Graduate studies in
  mathematics, 38, American mathematical society, Providence, RI, 2002.

\bibitem{Ik-Wat}
N.~Ikeda and S.~Watanabe, \emph{Stochastic differential equations and diffusion
  processes}, second ed., North-Holland Mathematical Library, 24, North-Holland
  Publishing Co., Amsterdam-New York; Kodansha, Ltd., Tokyo, 1989.

\bibitem{Kend}
W.~Kendall, \emph{Nonnegative {R}icci curvature and the {B}rownian coupling
  property}, Stochastics \textbf{19} (1986), 111--129.

\bibitem{Kend_survey}
W.S. Kendall, \emph{From stochastic parallel transport to harmonic maps}, New
  directions in Dirichlet forms, AMS/IP Studies in Advanced Mathematics, 8,
  Amer. Math. Soc., Providence, RI; International Press, Cambridge, MA, 1998,
  pp.~49--115.

\bibitem{Kurtz_appr-semigr}
T.G. Kurtz, \emph{Extensions of {T}rotter's operator semigroup approximation
  theorems}, J. Funct. Anal. \textbf{3} (1969), 111--132.

\bibitem{K8}
K.~Kuwada, \emph{Couplings of the {B}rownian motion via discrete apporoximation
  under lower {R}icci curvature bounds}, Probabilistic Approach to Geometry
  (Tokyo), Adv. Stud. Pure Math. 57, Math. Soc. Japan, 2010, pp.~273--292.

\bibitem{K-Phili2}
K.~Kuwada and R.~Philipowski, \emph{Coupling of {B}rownian motion and
  {P}erelman's {$\mathcal{L}$}-functional}, In preparation.

\bibitem{K-Phili}
\bysame, \emph{Non-explosion of diffusion processes on manifolds with
  time-dependent metric}, To appear in Math. Z.

\bibitem{McC-Topp_Wass-RF}
R.J. McCann and P.~Topping, \emph{Ricci flow, entropy and optimal
  transportation}, Amer. J. Math. \textbf{132} (2010), 711--730.

\bibitem{Oshima_tdepDF}
Y.~Oshima, \emph{Time-dependent {D}irichlet forms and related stochastic
  calculus}, Infin. Dimens. Anal. Quantum Probab. Relat. Top. \textbf{7}
  (2004), no.~2, 281--316.

\bibitem{Phili_sem}
R.~Philipowski, \emph{Coupling of diffusions on manifolds with 
  time-dependent metric}, Seminar talk at Universtit{\"a}t Bonn (2009). 

\bibitem{Qian:BonnetMyers}
Z.-M. Qian, \emph{Estimates for weighted volumes and applications}, Quart. J.
  Math. Oxford Ser. (2) \textbf{48} (1997), 235--242.

\bibitem{Str-Var}
D.W. Stroock and S.R.S. Varadhan, \emph{Multidimensional diffusion processes},
  Grundlehren der Mathematischen Wissenschaften, 233, Springer-Verlag, Berlin
  and New York, 1979.

\bibitem{Topp_Lopt}
P.~Topping, \emph{{$\mathcal{L}$}-optimal transportation for {R}icci flow}, J.
  Reine Angew. Math. \textbf{636} (2009), 93--122.

\bibitem{Renes_poly}
M.-K. von Renesse, \emph{Intrinsic coupling on {R}iemannian manifolds and
  polyhedra}, Electron. J. Probab. \textbf{9} (2004), no.~14, 411--435.

\bibitem{Wang94}
F.-Y. Wang, \emph{Successful couplings of nondegenerate diffusion processes on
  compact manifolds}, Acta. Math. Sinica. \textbf{37} (1994), no.~1, 116--121.

\bibitem{Wang_book05}
\bysame, \emph{Functional inequalities, {M}arkov semigroups, and spectral
  theory}, Mathematics Monograph Series 4, Science Press, Beijing, China, 2005.

\bibitem{Will-Wong}
D.~Willett and J.S.W. Wong, \emph{On the discrete analogues of some
  generalizations of {G}ronwall's inequality}, Monatsh. Math. \textbf{69}
  (1965), 362--367.

\bibitem{ZhangQS_HK-RF-Poinc}
Qi~S. Zhang, \emph{Heat kernel bounds, ancient {$\kappa$} soliton and the
  {P}oincar{\'{e}} conjecture}, J. Funct. Anal. \textbf{258} (2010), no.~4,
  1225--1246.

\end{thebibliography}
\providecommand{\bysame}{\leavevmode\hbox to3em{\hrulefill}\thinspace}
\providecommand{\MR}{\relax\ifhmode\unskip\space\fi MR }
\providecommand{\MRhref}[2]{%
  \href{http://www.ams.org/mathscinet-getitem?mr=#1}{#2}
}
\providecommand{\href}[2]{#2}

\vspace{1cm}
\begin{flushright}
Kazumasa Kuwada\\
\vspace{.3cm}
\small
Graduate School of Humanities and Sciences 
\\
Ochanomizu University 
\\
Ohtsuka 2-1-1, Bunkyo-ku, Tokyo 112-8610, Japan
\\
\vspace{.2cm}
\textit{e-mail}: \texttt{kuwada.kazumasa@ocha.ac.jp}
\\
\textit{URL}: \texttt{http://www.math.ocha.ac.jp/kuwada}
\vspace{1em}
\\
\normalsize
\textbf{Present address:} 
\\
\small 
Institut f\"ur Angewandte Mathematik 
\\
Universit\"at Bonn 
\\
Endenicher Allee 60, 53115 Bonn, Germany 
\\
\vspace{.2cm}
\textit{email}: \texttt{kuwada@iam.uni-bonn.de}
\end{flushright}

\end{document}